\documentclass[english,final]{siamart190516}
\usepackage[T1]{fontenc}
\usepackage[latin9]{inputenc}
\usepackage{refstyle}
\usepackage{mathrsfs}
\usepackage{amstext}
\usepackage{amssymb}
\usepackage{graphicx}

\makeatletter

\AtBeginDocument{\providecommand\thmref[1]{\ref{thm:#1}}}
\AtBeginDocument{\providecommand\exaref[1]{\ref{exa:#1}}}
\AtBeginDocument{\providecommand\figref[1]{\ref{fig:#1}}}
\AtBeginDocument{\providecommand\secref[1]{\ref{sec:#1}}}
\AtBeginDocument{\providecommand\defref[1]{\ref{def:#1}}}
\AtBeginDocument{\providecommand\corref[1]{\ref{cor:#1}}}
\AtBeginDocument{\providecommand\lemref[1]{\ref{lem:#1}}}
\AtBeginDocument{\providecommand\factref[1]{\ref{fact:#1}}}
\AtBeginDocument{\providecommand\appref[1]{\ref{app:#1}}}
\RS@ifundefined{subsecref}
  {\newref{subsec}{name = \RSsectxt}}
  {}
\RS@ifundefined{thmref}
  {\def\RSthmtxt{theorem~}\newref{thm}{name = \RSthmtxt}}
  {}
\RS@ifundefined{lemref}
  {\def\RSlemtxt{lemma~}\newref{lem}{name = \RSlemtxt}}
  {}

\numberwithin{equation}{section}

\usepackage{url}
\usepackage{xcolor}
\newref{prop}{name=Proposition~,Name=Proposition~}
\newref{thm}{name=Theorem~,Name=Theorem~}
\newref{lem}{name=Lemma~,Name=Lemma~}
\newref{def}{name=Definition~,Name=Definition~}
\newref{cor}{name=Corollary~,Name=Corollary~}
\newsiamthm{fact}{Fact}
\newref{fact}{name=Fact~,Name=Fact~}
\newref{sec}{name=Section~,Name=Section~}
\newref{app}{name=Appendix~,Name=Appendix~}
\newsiamthm{example}{Example}
\newref{exa}{name=Example~,Name=Example~}
\newref{fig}{name=Fig.~,Name=Fig.~}

\newref{asm}{name=Assumption~,Name=Assumption~}

\makeatother

\usepackage{babel}

\begin{document}
\title{Sharp Global Guarantees for Nonconvex Low-rank Recovery in the Noisy
Overparameterized Regime\thanks{Financial support for this work was provided by NSF CAREER Award ECCS-2047462
and ONR Award N00014-24-1-2671.}}
\author{Richard Y. Zhang\thanks{Dept. of Electrical and Computer Engineering, University of Illinois
at Urbana-Champaign, 306 N Wright St, Urbana, IL 61801, \url{ryz@illinois.edu}}}
\maketitle
\begin{abstract}
Recent work established that rank overparameterization eliminates
spurious local minima in nonconvex low-rank matrix recovery under
the restricted isometry property (RIP). But this does not fully explain
the practical success of overparameterization, because real algorithms
can still become trapped at nonstrict saddle points (approximate second-order
points with arbitrarily small negative curvature) even when all local
minima are global. Moreover, the result does not accommodate for noisy
measurements, but it is unclear whether such an extension is even
possible, in view of the many discontinuous and unintuitive behaviors
already known for the overparameterized regime. In this paper, we
introduce a novel proof technique that unifies, simplifies, and strengthens
two previously competing approaches---one based on escape directions
and the other based on the inexistence of counterexample---to provide
sharp global guarantees in the noisy overparameterized regime. We
show, once local minima have been converted into global minima through
slight overparameterization, that near-second-order points achieve
the same minimax-optimal recovery bounds (up to small constant factors)
as significantly more expensive convex approaches. Our results are
sharp with respect to the noise level and the solution accuracy, and
hold for both the symmetric parameterization $XX^{T}$, as well as
the asymmetric parameterization $UV^{T}$ under a balancing regularizer;
we demonstrate that the balancing regularizer is indeed necessary. 
\end{abstract}

\global\long\def\R{\mathbb{R}}%
\global\long\def\S{\mathbb{S}}%
\global\long\def\rank{\mathrm{rank}}%
\global\long\def\nnz{\mathrm{nnz}}%
\global\long\def\e{\mathbf{e}}%
\global\long\def\J{\mathbf{J}}%
\global\long\def\A{\mathbf{A}}%
\global\long\def\G{\mathbf{G}}%
\global\long\def\H{\mathbf{H}}%
\global\long\def\vect{\mathrm{vec}}%
\global\long\def\mat{\mathrm{mat}}%
\global\long\def\ub{\mathrm{ub}}%
\global\long\def\lb{\mathrm{lb}}%
\global\long\def\tr{\mathrm{tr}}%
\global\long\def\eqdef{\overset{\mathrm{def}}{=}}%
 
\global\long\def\T{\mathrm{T}}%
\global\long\def\N{\mathrm{N}}%
\global\long\def\one{\mathbf{1}}%
\global\long\def\half{{\textstyle \frac{1}{2}}}%
\global\long\def\f{\mathbf{f}}%

\global\long\def\colsp{\mathrm{colsp}}%
\global\long\def\op{\mathrm{op}}%
\global\long\def\nuc{\mathrm{nuc}}%

\global\long\def\Prim{\operatorname{Prim}}%
\global\long\def\Dual{\operatorname{Dual}}%

\global\long\def\rip{\operatorname{RIP}}%
\global\long\def\ripd{\rip(\delta,r+r^{\star})}%

\global\long\def\rwc{\operatorname{RWC}}%
\global\long\def\rwcd{\rwc(\delta,r+r^{\star})}%

\global\long\def\grad{\operatorname{grad}}%
\global\long\def\Hess{\operatorname{Hess}}%

\global\long\def\orth{\operatorname{orth}}%
\global\long\def\cond{\operatorname{cond}}%

\global\long\def\diag{\operatorname{diag}}%

\global\long\def\Diag{\operatorname{Diag}}%

\global\long\def\offdiag{\operatorname{offdiag}}%

\global\long\def\inner#1#2{\left\langle #1,#2\right\rangle }%
\global\long\def\Id{\mathrm{Id}}%
\global\long\def\cF{\mathscr{F}}%
\global\long\def\cG{\mathscr{G}}%
\global\long\def\cH{\mathscr{H}}%

\global\long\def\hD{\hat{D}}%
\global\long\def\hX{\hat{X}}%
\global\long\def\hZ{\hat{Z}}%
\global\long\def\hU{\hat{U}}%

\global\long\def\cM{\mathcal{M}}%
\global\long\def\AA{\mathcal{A}}%
\global\long\def\BB{\mathcal{B}}%
\global\long\def\II{\mathcal{I}}%
\global\long\def\HH{\mathcal{H}}%
\global\long\def\LL{\mathcal{L}}%
\global\long\def\UU{\mathcal{U}}%

\section{Introduction}

Low-rank matrix recovery seeks to estimate an unknown $n_{1}\times n_{2}$
matrix $M^{\star}$ of low-rank $r^{\star}$ from noisy measurements
$b\approx\AA(M^{\star})$ made by a known linear operator $\AA$.
Convex methods achieve minimax-optimal recovery~\cite{recht2010guaranteed,candes2011tight}
but are often too computationally expensive for real-world data. Instead,
it is more common to apply a cheap gradient-based algorithm to the
nonconvex least-squares problem
\begin{equation}
\text{minimize}\quad\|\AA(UV^{T})-b\|^{2}\quad\text{over }U\in\R^{n_{1}\times r},V\in\R^{n_{2}\times r}.\label{eq:prob1}
\end{equation}
The key feature of this nonconvex approach is that it reduces the
number of parameters from quadratic $n_{1}n_{2}$ down to linear $(n_{1}+n_{2})r$
where $r\ge r^{\star}$ is the model rank. 

While setting $r=r^{\star}$ would indeed minimize the total number
of parameters, practitioners often \emph{overparameterize} the model
rank $r>r^{\star}$, as numerical evidence suggests that doing so
produces a more benign optimization landscape and reduces the risk
of getting stuck at spurious local minima~\cite{rosen2019se,chen2020accelerating,zhang2023preconditioned,chiu2023tight}.
Recently, \cite{zhang2022improved} made progress towards a theoretical
explanation of this phenomenon under the \emph{restricted isometry
property} (RIP).   In what follows, we write $\inner EF=\tr(E^{T}F)$
and $\|E\|=\sqrt{\inner EE}$ to denote the matrix Euclidean (i.e.
Frobenius) inner product and norm.
\begin{definition}[RIP]
Denote $\rip(\delta,k)$ as the set of all linear maps $\AA$ satisfying
the $(\delta,k)$\nobreakdash-\emph{restricted isometry property}
for $0\le\delta<1$ and $k\ge1$: 
\[
\rank(E)\le k\quad\implies\quad(1-\delta)\|E\|^{2}\le\|\AA(E)\|^{2}\le(1+\delta)\|E\|^{2}.
\]
\end{definition}

When RIP holds with sufficiently small $\delta\approx0$, the symmetric
formulation of (\ref{eq:prob1}) with $U=V$ is well-known to exhibit
a \emph{benign landscape}: all second-order points (and thus all local
minima) are global minima that recover the ground truth~\cite{bhojanapalli2016global,ge2017nospurious}.
This follows from a perturbative analysis:
\[
\|\AA(XX^{T})-b\|^{2}\quad\overset{\text{noisy meas}}{\approx}\quad\|\AA(XX^{T}-M^{\star})\|^{2}\quad\overset{\text{RIP}}{\approx}\quad\|XX^{T}-M^{\star}\|^{2},
\]
where the final expression, the squared error norm, can be verified
to have a benign landscape~\cite{ge2017nospurious}. But if $\delta$
is too large, then it becomes possible for spurious local minima to
emerge~\cite{zhang2018much,zhang2019sharp}. Instead, \cite{zhang2022improved}
showed that overparameterizing the model rank $r>r^{\star}$ eliminates
these spurious local minima, hence improving the landscape and eventually
rendering it benign at the threshold $r/r^{\star}>[\delta/(1-\delta)]^{2}$. 
\begin{theorem}[{\cite[Corollary~1.5]{zhang2022improved}}]
\label{thm:zhang}Let $M^{\star}\in\R^{n\times n}$ satisfy $M^{\star}\succeq0$
and $\rank(M^{\star})\le r^{\star}$, and let $\AA\in\rip(\delta,k)$.
For $r$ satisfying $r^{\star}\le r<n$, define $f:\R^{n\times r}\to\R$
such that 
\[
f(X)=\half\|\AA(XX^{T})-b\|^{2}\quad\text{where }b=\AA(M^{\star}).
\]
If $r/r^{\star}>[\delta/(1-\delta)]^{2}$ and $k\ge r+r^{\star}$,
then every exact second-order point exactly recovers the ground truth:
\[
\nabla f(X)=0,\quad\nabla^{2}f(X)\succeq0\quad\iff\quad XX^{T}=M^{\star}.
\]
If $r/r^{\star}\le[\delta/(1-\delta)]^{2}$, then for every $k\ge1$,
there exists a counterexample that admits a spurious second-order
point $X_{0}$ with error $\|X_{0}X_{0}^{T}-M^{\star}\|>\|M^{\star}\|$.
\end{theorem}

While $f$ may still have a benign landscape below the threshold $r/r^{\star}\le[\delta/(1-\delta)]^{2}$,
no RIP-based guarantee is possible in that regime, because one cannot
distinguish $f$ from the counterexample with spurious local minima.
Indeed, the elimination of such counterexamples beyond the threshold
serves as strong evidence of overparameterization\textquoteright s
ability to improve the landscape. However, \thmref{zhang} still does
not fully explain the practical and algorithmic success of overparameterization.
Three critical gaps arise: (i) local minima vs approximate second-order
points; (ii) noiseless vs noisy measurements; (iii) symmetric vs asymmetric
parameterizations. 

First, real algorithms cannot compute local minima or even exact second-order
points, but only \emph{approximate second-order points} that satisfy
the second-order optimality conditions within some small tolerance.
This makes them susceptible to stalling at a spurious point with arbitrarily
small negative curvature, called a \emph{nonstrict saddle point},
even when the landscape is benign. The following is the $r=r^{\star}=1$
and $n=2$ instance of \exaref{sym} later in the paper; it also generalizes
\cite{zhang2018much}. 
\begin{figure}
\begin{minipage}[b][1\totalheight][t]{0.5\columnwidth}%
\includegraphics[width=1\columnwidth]{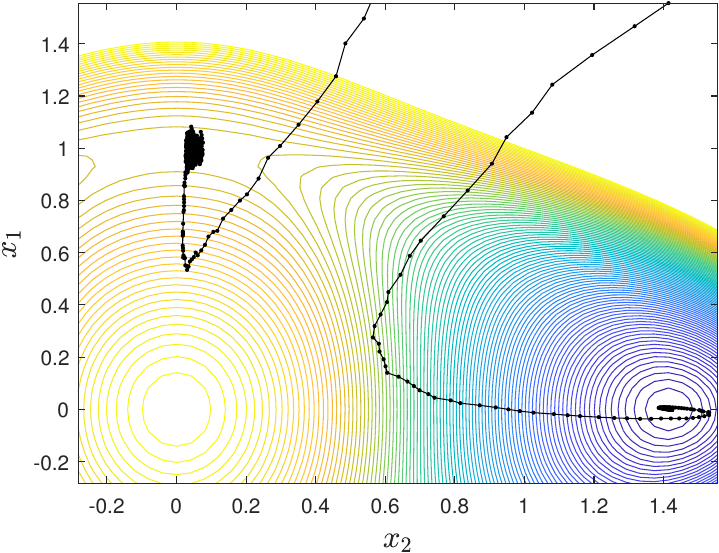}%
\end{minipage}%
\begin{minipage}[b][1\totalheight][t]{0.5\columnwidth}%
\includegraphics[width=0.98\columnwidth]{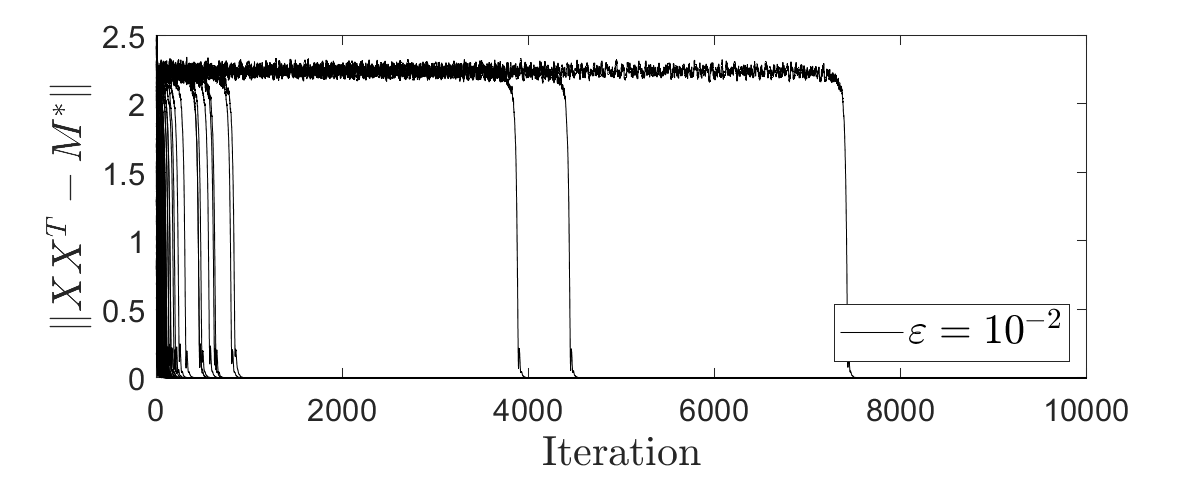}

\includegraphics[width=0.98\columnwidth]{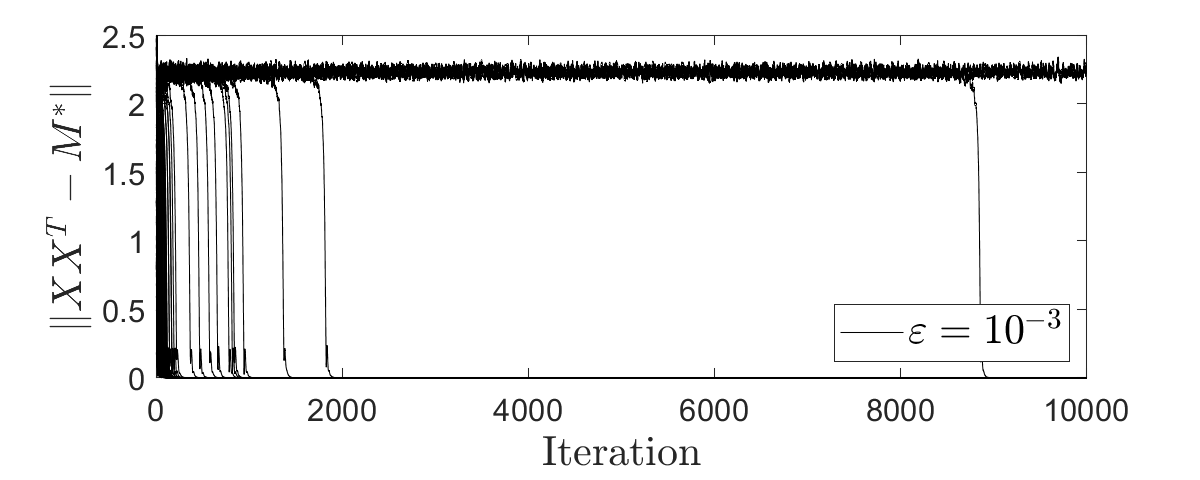}%
\end{minipage}

\caption{\textbf{Nonstrict saddle point can stall SGD even when the landscape
is benign.} \emph{Left.} For \exaref{approx} with $\varepsilon=10^{-2}$,
SGD with Nesterov momentum can stall at the nonstrict saddle point
point at $(1,0)$ after $10^{3}$ steps, even though the only local
minima (and thus global minima) lie at $(0,\pm\sqrt{2+\varepsilon})$.
\emph{Right top:} Extending to $10^{4}$ steps allows all 100 trials
to escape to the global minimum. \emph{Right bottom:} Lowering $\varepsilon=10^{-3}$
causes 2 out of 100 trials to stall after $10^{4}$ steps. (Experiment
details: Set $f(X)=\sum_{i=1}^{4}f_{i}(X)$ with $f_{i}(X)=\frac{1}{2}\left|\langle A_{i},XX^{T}-M^{\star}\rangle\right|^{2}$,
initialize $X\sim\mathcal{N}(0,I_{nr})$, $V=0$, update $V\gets\beta V-\alpha\nabla f_{i}(X)$
and $X\gets X+\beta V-\alpha\nabla f_{i}(X)$ using $\alpha=10^{-2}$,
$\beta=0.9$, increment index $i$ modulo 4, and shuffle every epoch$=4$
iterations.) \protect\label{fig:example}}
\end{figure}

\begin{example}[Failure by nonstrict saddle point]
\label{exa:approx}For $\varepsilon>0,$ define $f(x)=\frac{1}{2}\|\AA(xx^{T}-M^{\star})\|^{2}$
where $M^{\star}$ and $\AA(M)=[\inner{A_{i}}M]_{i=1}^{4}$ are given
\begin{gather*}
M^{\star}=\begin{bmatrix}0 & 0\\
0 & 2+\varepsilon
\end{bmatrix},\quad A_{1}=\sqrt{\frac{1+\delta}{2}}\begin{bmatrix}1 & 0\\
0 & 1
\end{bmatrix},\quad A_{2}=\sqrt{\frac{1-\delta}{2}}\begin{bmatrix}1 & 0\\
0 & -1
\end{bmatrix},\\
A_{3}=\sqrt{1+\delta}\begin{bmatrix}0 & 0\\
1 & 0
\end{bmatrix},\quad A_{4}=\sqrt{1+\delta}\begin{bmatrix}0 & 1\\
0 & 0
\end{bmatrix},\quad\delta=\frac{1}{2+\varepsilon}.
\end{gather*}
It is easy to verify that $\AA\in\rip(\delta,2)$ with $\delta<\frac{1}{2}$,
so according to \thmref{zhang} (and also \cite[Theorem~3]{zhang2019sharp}),
all exact second-order points (and hence all local minima) are global
minima that exactly recover the ground truth:
\[
\nabla f(x)=0,\quad\nabla^{2}f(x)\succeq0\quad\iff\quad xx^{T}=M^{\star}.
\]
But $x_{0}=[1;0]$ is nonstrict saddle point, i.e. an approximate
second-order point that fails to achieve near-recovery: 
\[
\nabla f(x_{0})=\begin{bmatrix}0\\
0
\end{bmatrix},\quad\nabla^{2}f(x_{0})=\begin{bmatrix}4 & 0\\
0 & -2\varepsilon(1+\frac{1}{2+\varepsilon})
\end{bmatrix},\quad\|x_{0}x_{0}^{T}-M^{\star}\|>\|M^{\star}\|.
\]
\end{example}

\figref{example} illustrates how SGD could fail on \exaref{approx}
by stalling at a nonstrict saddle point. This example highlights an
important point: while \thmref{zhang} guarantees the inexistence
of spurious local minima, real algorithms further require the inexistence
of nonstrict saddle points to reliably recover the ground truth. In
order to derive rigorous convergence rates or prove complexity bounds
based on rank and RIP, one must further show that all approximate
second-order points yield near-recovery.

Second, an exact recovery guarantee like \thmref{zhang} cannot accommodate
\emph{noisy measurements}, which are inevitable due to physical limitations
and numerical round-off, but it is unclear whether the extension to
inexact recovery is even possible. In the overparameterized regime
$r>r^{\star}$, the iterates of a recovery algorithm must approach
the boundary of the Riemannian manifold of rank-$r$ matrices, where
the local curvature grows to be infinite in the limit. Most existing
inexact recovery guarantees~\cite{bhojanapalli2016dropping,bhojanapalli2016global,ge2017nospurious,zheng2015convergent,tu2016low}
rely on a finite local curvature to control the recovery error, and
so become vacuous once $r>r^{\star}$. Indeed, the infinite curvature
at the boundary causes many discontinuous and unintuitive behavior
to manifest, such as the exponential slowdown of gradient algorithms~\cite{zhuo2024computational,xiong2024over},
and the appearance of ``apocalypses'' that cause first-order algorithms
to fail~\cite{levin2023finding}. Even though minimax-optimal recovery
is achieved in the overparameterized regime by specific nonconvex
algorithms~\cite{stoger2021small,zhang2021preconditioned,xu2023power,zhuo2024computational},
it is still conceivable that a worst-case noise perturbation to a
vulnerable algorithm could cause recovery error to blow up exponentially,
such that the benefits of overparameterization would no longer materialize
across all algorithms. 

Third, as an important practical point, \thmref{zhang} does not cover
the \emph{asymmetric parameterization} $UV^{T}$ that is more commonly
used to recover general low-rank matrices that are possibly indefinite,
nonsymmetric, or even nonsquare. In fact, the asymmetric parameterization
$UV^{T}$ can be more preferable even when the ground truth $M^{\star}$
is known to be symmetric positive semidefinite~\cite{xiong2024over}.

Unfortunately, bridging these three gaps has proved unexpectedly difficult,
owing to a sharpness--generality trade-off between existing proof
techniques. In one direction, global guarantees like \thmref{zhang}
have long been generalized to approximate second-order points~\cite{ge2017nospurious,zhu2018global,li2019non},
noisy measurements~\cite{bhojanapalli2016global,ge2017nospurious},
and the asymmetric parameterization~\cite{park2017non,ge2017nospurious},
but none of these have been able to demonstrate an improvement with
overparameterization. The fundamental barrier is the conservatism
of the underlying proof technique, based on the \emph{existence of
an escape direction}~\cite{zheng2015convergent,tu2016low,zhu2018global}
(also known as a \emph{direction of improvement}~\cite{ge2017nospurious}),
as it is not sharp enough to capture a dependence on overparameterization.
\thmref{zhang} was the first to uncover an unambiguous improvement
with overparameterization, precisely because it had fully sharpened
the conservatism of prior work using a novel proof technique based
on the \emph{inexistence of counterexamples}~\cite{zhang2019sharp}.
But as we explain in detail in \secref{duality}, the sharpness of
this proof technique also makes it very fragile and difficult to generalize
to broader settings. 

\subsection{Main results}

This paper presents a simplified proof technique that eliminates the
sharpness--generality tradeoff across the three gaps in the existing
literature. Our critical insight is that the two previous competing
approaches---one based on the existence of an escape direction and
the other based on the inexistence of a counterexample---are in fact
\emph{strong Lagrangian duals} of each other. While the existence
of an escape direction obviously implies the inexistence of a counterexample,
strong Lagrangian duality further ensures that an escape direction
exists if and only if a counterexample does not exist (\thmref{strongdual}).
Critically, the sharpness of \thmref{zhang}, proved using the inexistence
of counterexamples, implies via strong duality the existence of corresponding
\emph{sharp escape directions}, that will fully sharpen existing proofs
based on escape directions, that in turn easily generalize across
the gaps. 

Our main technical contribution is to explicitly identify these sharp
escape directions (\defref{directions}), and to use these to provide
a natural generalization of \thmref{zhang} to noisy measurements
and approximate second-order points. In fact, the generalization is
sharp with respect to the noise level $\epsilon_{0}$ and the accuracy
parameters $\epsilon_{1},\epsilon_{2}$, up to small absolute constant
multipliers. Below, $\|\cdot\|_{\op}$ denotes the matrix operator
norm (i.e. the spectral norm). 
\begin{theorem}[Symmetric parameterization]
\label{thm:sym}Let $M^{\star}\in\R^{n\times n}$ satisfy $M^{\star}\succeq0$
and $\rank(M^{\star})\le r^{\star}$, and let $\AA\in\rip(\delta,k)$.
For $r$ satisfying $r^{\star}\le r<n$, define $f:\R^{n\times r}\to\R$
such that 
\[
f(X)=\half\|\AA(XX^{T})-b\|^{2}\quad\text{where }b=\AA(M^{\star})+w,\;\|\AA^{T}(w)\|_{\op}\le\epsilon_{0},
\]
and let $X\in\R^{n\times r}$ denote an $(\epsilon_{1},\epsilon_{2})$-approximate
second-order point
\[
\inner{\nabla f(X)}D\ge-\epsilon_{1}\|DX^{T}\|,\quad\inner{\nabla^{2}f(X)[D]}D\ge-\epsilon_{2}\|DX^{T}\|^{2}\quad\text{for all }D\in\R^{n\times r}.
\]
If $r/r^{\star}>[(1+\half\epsilon_{2})\delta/(1-\delta)]^{2}$ and
$k\ge r+r^{\star}$, then $X$ nearly recovers the ground truth:
\begin{gather*}
\|XX^{T}-M^{\star}\|\le\left(\frac{\epsilon_{1}}{2}+\epsilon_{0}\sqrt{r+r^{\star}}\right)\cdot\left(\frac{1}{1+(1+\half\epsilon_{2})\sqrt{r^{\star}/r}}-\delta\right)^{-1}.
\end{gather*}
If $r/r^{\star}\le[(1+\frac{1}{4}\epsilon_{2})\delta/(1-\delta)]^{2}$,
then for every $k\ge1$, there exists a counterexample that admits
a spurious $(\epsilon_{1},\epsilon_{2})$-approximate second-order
point $X_{0}$ with error $\|X_{0}X_{0}^{T}-M^{\star}\|>\|M^{\star}\|$.
\end{theorem}

The noise model $\|\AA^{T}(w)\|_{\op}\le\epsilon_{0}$ and the local
norm $\|DX^{T}\|$ used to define approximate second-order points
in \thmref{sym} are chosen specifically to yield a sharp dependence.
 There exist algorithms that directly compute approximate second-order
points in the local norm~\cite{uschmajew2020critical,zhang2023preconditioned,xu2023power},
which is closely related to the Euclidean metric on the Riemannian
manifold of rank-$r$ positive semidefinite matrices; see~\cite{mishra2012riemannian,uschmajew2020critical}
and also \cite[Section 7.5]{boumal2023introduction}. By giving up
the sharp dependence on $\epsilon_{1}$ and $\epsilon_{2}$, \thmref{sym}
can be translated into the standard Euclidean norm.
\begin{corollary}
\label{cor:sym}Under the same condition as \thmref{sym}, let $w\sim\mathcal{N}(0,\sigma^{2}\cdot I_{m})$.
If $r/r^{\star}>[\delta/(1-\delta)_{+}]^{2}$ and $k\ge r+r^{\star}$,
then any approximate second-order point $X$ that satisfies $\|\nabla f(X)\|\le\epsilon_{1}$
and $\nabla^{2}f(X)\succeq-\epsilon_{2}I$ has error
\begin{gather*}
\|XX^{T}-M^{\star}\|\le\frac{15\sigma\sqrt{n(r+r^{\star})}}{\Delta}+\sqrt{\frac{\epsilon_{1}R+\epsilon_{2}R^{2}}{\Delta}},
\end{gather*}
where $\Delta=\frac{1}{1+\sqrt{r^{\star}/r}}-\delta$ and $R=\sqrt{(1+\delta)(1+\sqrt{r/r^{\star}})}\max\{\|X\|,\|Z\|\}$,
with probability at least $1-12^{-2n+1}$ over the randomness of the
noise vector $w$. 
\end{corollary}

\corref{sym} allows us to easily derive the algorithmic implications
of rank overparameterization, including convergence rates and complexity
bounds based on rank and RIP parameters. Indeed, numerous algorithms
are capable of computing the specified approximate second-order point
$X$: stochastic gradient descent (SGD) with suitable perturbations
requires at most $(\epsilon_{1}^{-4}+\epsilon_{2}^{-8})\cdot\mathrm{polylog}(\epsilon_{1}^{-1},\epsilon_{2}^{-1},n,r)$
iterations~\cite{jin2021nonconvex}, while cubically-regularized
Newton requires at most $O(\epsilon_{1}^{-3/2}+\epsilon_{2}^{-3})$
iterations~\cite{nesterov2006cubic}. Under standard RIP assumptions
and adopting a Gaussian noise model, \corref{sym} says that slightly
overparameterizing $r/r^{\star}>[\delta/(1-\delta)]^{2}$ is both
sufficient and necessary for minimax-optimal recovery $XX^{T}\approx M^{\star}$.
For example, choosing a model rank of $r=[\delta/(1-\delta)]^{2}r^{\star}+1$,
it takes perturbed SGD at most $\tilde{O}(\sigma^{-16})$ iterations
to arrive at an estimate $X$ with error $\|XX^{T}-M^{\star}\|\le\frac{20}{\Delta^{2}}\sigma\sqrt{nr^{\star}}$,
which is minimax optimal within a constant factor of $40/\Delta^{2}$~\cite[Theorem~2.6]{candes2011tight}.
(Here, we matched $\sqrt{\epsilon_{1}}=\sqrt{\epsilon_{2}}=\sigma$
and used $\sqrt{1+[\delta/(1-\delta)]^{2}}\le1+\delta/(1-\delta)=1/(1-\delta)\le1/\Delta$.)

Second, we show that the benefits of overparameterization for the
symmetric case $XX^{T}$ do indeed extend to the asymmetric case $UV^{T}$,
but only after the problem is augmented with a balancing regularizer.
The sufficient conditions extend almost verbatim from \thmref{sym}
and \corref{sym}, except that the RIP constant $\delta$ and the
noise level $\epsilon_{0}$ are effectively doubled. It turns out
that the counterexample used to prove the necessary condition also
extends to the asymmetric case.
\begin{theorem}[Asymmetric parameterization]
\label{thm:asym}Let $M^{\star}\in\R^{n_{1}\times n_{2}}$ satisfy
$\rank(M^{\star})\le r^{\star}$, and let $\AA\in\rip(\delta,k)$.
For $r\ge r^{\star}$, define $g:\R^{n\times r}\to\R$ where $n=n_{1}+n_{2}$
such that 
\[
g([U;V])=2\|\AA(UV^{T})-b\|^{2}+\half\|U^{T}U-V^{T}V\|^{2}\;\text{where }b=\AA(M^{\star})+w,\;\|\AA^{T}(w)\|_{\op}\le\epsilon_{0},
\]
and let $X=[U;V]\in\R^{n\times r}$ denote an $(\epsilon_{1},\epsilon_{2})$-approximate
second-order point
\[
\inner{\nabla g(X)}D\ge-\epsilon_{1}\|DX^{T}\|,\quad\inner{\nabla^{2}g(X)[D]}D\ge-\epsilon_{2}\|DX^{T}\|^{2}\quad\text{for all }D\in\R^{n\times r}.
\]
If $r/r^{\star}>[2(1+\frac{1}{2}\epsilon_{2})\delta/(1-2\delta)_{+}]^{2}$
and $r+r^{\star}\le k$, then $X$ nearly recovers the ground truth:
\begin{gather*}
\|UV^{T}-M^{\star}\|\le\left(\frac{\epsilon_{1}}{2}+2\epsilon_{0}\sqrt{r+r^{\star}}\right)\cdot\left(\frac{1}{1+(1+\half\epsilon_{2})\sqrt{r^{\star}/r}}-2\delta\right)^{-1}
\end{gather*}
If $r/r^{\star}\le[(1+\frac{1}{4}\epsilon_{2})\delta/(1-\delta)]^{2}$,
then for every $k\ge1$ and $t\ge0$, there exists a counterexample
that admits a spurious $(\epsilon_{1},\epsilon_{2})$-approximate
second-order point $X_{0}=[U_{0};V_{0}]$ with error $\|U_{0}V_{0}^{T}-M^{\star}\|>\|M^{\star}\|$.
\end{theorem}

\begin{corollary}
\label{cor:nonsym}Under the same conditions as \thmref{asym}, suppose
that $w\sim\mathcal{N}(0,\sigma^{2}\cdot I_{m})$. If $r/r^{\star}>[2\delta/(1-2\delta)_{+}]^{2}$
and $k\ge r+r^{\star}$, then any approximate second-order point $X=[U;V]$
that satisfies $\|\nabla f(X)\|\le\epsilon_{1}$ and $\nabla^{2}f(X)\succeq-\epsilon_{2}I$
has error
\begin{gather*}
\|UV^{T}-M^{\star}\|\le\frac{20\sigma\sqrt{n(r+r^{\star})}}{\Delta}+\sqrt{\frac{\epsilon_{1}R+\epsilon_{2}R^{2}}{\Delta}},
\end{gather*}
where $\Delta=\frac{1}{1+\sqrt{r^{\star}/r}}-2\delta$ and $R=\sqrt{(1+2\delta)(1+\sqrt{r/r^{\star}})}\max\{\|X\|,\|Z\|\}$,
with probability at least $1-12^{-n_{1}-n_{2}+1}$ over the randomness
of the noise vector $w$. 
\end{corollary}

Analogous to the symmetric case before, \corref{nonsym} says overparameterizing
$r/r^{\star}>[2\delta/(1-2\delta)]^{2}$ is sufficient for near-second-order
points to achieve minimax-optimal recovery $XX^{T}\approx M^{\star}$,
while $r/r^{\star}>[\delta/(1-\delta)]^{2}$ is necessary for such
a recovery guarantee to be possible at all. Note that augmenting with
the balancing regularizer $\|U^{T}U-V^{T}V\|^{2}$ does not change
the nature of the problem, because every low-rank matrix $M^{\star}$
admits a balanced factorization $M^{\star}=U^{\star}V^{\star T}$
with $U^{\star T}U^{\star}=V^{\star T}V^{\star}$. On the other hand,
\thmref{asym} and \corref{nonsym} would not have been possible without
this regularizer term, due to the following counterexample. 
\begin{example}[Necessity of balancing regularizer]
\label{exa:bal}Let $M^{\star}=e_{1}e_{1}^{T}$ and $\AA\in\rip(0,n)$.
For $r$ satisfying $1\le r<n$, define $f:\R^{2n\times r}\to\R$
such that 
\[
f([U;V])=\|\AA(UV^{T})-b\|^{2}\quad\text{where }b=\AA(M^{\star}).
\]
Then, $X=[U;V]$ with $U=(1/\sqrt{\epsilon})[e_{2},e_{3},\dots,e_{r+1}]$
and $V=0$ is a spurious approximate second-order point with $\nabla f(X)=0$
and $\nabla^{2}f(X)\succeq-\epsilon I$ but $\|UV^{T}-M^{\star}\|=1$.
\end{example}

Hence, the fact that existing algorithms have been able to achieve
similar behavior as \thmref{asym} without the regularizer indicate
that they enjoy some kind of implicit regularization that allow the
factors to remain balanced. This implicit balancing behavior has been
rigorously established for gradient descent starting from a (very)
small random initialization~\cite{du2018algorithmic,ye2021global,jiang2023algorithmic,soltanolkotabi2023implicit},
though it remains future work to understand why the behavior also
seems to manifest with any arbitrary initialization.

\subsection{Limitations and future work}

Having addressed the three main gaps highlighted in the introduction,
the remaining challenge lies in the inherent strength of the RIP assumption.
Although low-rank matrix recovery under RIP is standard in the literature~\cite{bhojanapalli2016global,ge2017nospurious,chi2019nonconvex},
a common critique is that it requires \emph{dense} measurements, in
the sense that $\AA(e_{i}e_{j}^{T})\ne0$ must hold for all $i,j$.
Consequently, our results do not directly apply to \emph{sparse} measurement
problems like matrix completion or phase retrieval. Nevertheless,
real-world applications with dense measurements suspected of satisfying
RIP do exist~\cite{lang2024interacting}, even though verifying RIP
formally remains NP-hard. In the literature, the RIP framework is
often viewed as a tractable starting point for developing proof techniques
that can later extend to non-RIP settings~\cite{ge2017nospurious,chi2019nonconvex},
rather than a literal reflection of real-world measurements.

A further limitation is the requirement $k\ge r+r^{\star}$ on the
RIP rank parameter, which forces the sample complexity $m$ to scale
with the model rank $r$ rather than the true rank $r^{\star}$. Since
the statistical problem itself depends only on $r^{\star}$, this
dependence on $r$ must be an artifact of the optimization approach.
This challenge is not unique to our work; all prior benign landscape
guarantees~\cite{bhojanapalli2016global,ge2017nospurious,zhu2018global,li2019non,mcrae2024low}
impose the same rank requirement. The fundamental issue is in approximating
$\|\mathcal{A}(XX^{T}-M^{\star})\|^{2}\approx\|XX^{T}-M^{\star}\|^{2}$
uniformly over all $X\in\mathbb{R}^{n\times r}$, as this holds only
under rank-$k$ RIP with $k\ge r+r^{\star}$. Intuitively, if there
exist multiple unrelated $XX^{T}$ yielding the same residual $\AA(XX^{T}-M^{\star})$,
then the loss landscape cannot be benign, and even global minimizers
may not recover the ground truth. We conjecture that $k\geq r+r^{\star}$
is both necessary and sufficient for a globally benign landscape;
a rigorous justification would require constructing counterexamples
that yield spurious second-order points when $k<r+r^{\star}$, in
line with \exaref{sym} and \exaref{asym} later in the paper. As
a partial (though not entirely satisfying) workaround, note that when
$r/r^{\star}=\Theta(1)$, both $r$ and $r^{\star}$ lie in the same
order, so $m$ still scales with $r^{\star}$ in practice.

Recently, recovery guarantees have been obtained under the weaker
condition $k\ge2r^{\star}$~\cite{stoger2021small,xu2023power,soltanolkotabi2023implicit,ma2023global,ma2024convergence},
but these rely on carefully chosen algorithmic trajectories and specific
initializations. The key idea is to start from a very small random
initial point, so that gradient descent exhibits \emph{incremental
learning}, increasing the rank of the overparameterized iterate $X$
one at a time until it reaches $r^{\star}$. Consequently, the approximation
$\|\AA(XX^{T}-M^{\star})\|^{2}\approx\|XX^{T}-M^{\star}\|^{2}$ only
requires rank-$2r^{\star}$ RIP, because $\rank(X)\le r^{\star}$
always holds despite overparameterization. This focus on a particular
trajectory contrasts with our primary goal of characterizing the entire
landscape. Our stronger rank-$(r+r^{\star})$ requirement highlights
the inherent challenge in maintaining a global benign landscape without
algorithmic restrictions. 

Therefore, an important direction for future work is to extend our
analysis to sparse measurements, and to rank-$2r^{\star}$ RIP. Based
on our results for the asymmetric parameterization, where a regularizer
was essential to render the landscape benign, we hypothesize that
similar regularization strategies could address these limitations.
 We speculate on the details of this at the end of the paper, in
\secref{conclusion}. 

\subsection{\protect\label{sec:related}Related and follow-up work}

Early theoretical work on nonconvex low-rank matrix recovery under
RIP (frequently called matrix sensing) focused on efficiently finding
a provably good initial guess for (\ref{eq:prob1}) using spectral
initialization~\cite{keshavan2010matrixa,jain2013low,chen2015solving,zheng2015convergent,tu2016low,sun2016guaranteed}.
Recently, a small random initialization was also found to be a provably
good initial guess~\cite{stoger2021small,xu2023power,soltanolkotabi2023implicit}.
In practice, however, any arbitrary initial guess seems to work just
as well as the provably good ones. In the RIP setting, \cite{bhojanapalli2016global}
was the first to prove a \emph{benign landscape} result, by showing
for $\AA\in\ripd$ with $\delta<1/5$ that all approximate second-order
points achieve near-recovery in the symmetric parameterization $XX^{T}$.
This was subsequently extended to the asymmetric parameterization
$UV^{T}$~\cite{park2017non,ge2017nospurious}, to non-RIP settings~\cite{ge2017nospurious},
and to general unconstrained low-rank optimization~\cite{zhu2018global,li2019non}.
But none of these prior results are able to explain why benign landscape
becomes more prevalent as the rank $r$ is overparameterized.

It was previously known that overparameterizing $r=\Omega(\sqrt{n})$
would cause functions like $f$ to generically have benign landscape~\cite{bhojanapalli2018smoothed,cifuentes2022polynomial},
but in practice, only $r=\Omega(r^{\star})$ seems to be needed for
the ground truth to be consistently and reliably recovered. To the
best of our knowledge, a weaker version of \thmref{zhang} appearing
in an earlier arXiv version~\cite{zhang2021sharp} of the present
paper was the first to establish benign landscape with slight overparameterization
$r=\Omega(r^{\star})$. Subsequently, \cite{zhang2022improved} strengthened
and generalized the result to unconstrained low-rank optimization,
albeit only for exact second-order points and noiseless measurements
and a symmetric parameterization. \cite{ma2023geometric} subsequently
used a continuity argument to derive a coarse error bound for noisy
measurements, but again only for exact second-order points. In contrast,
the present paper is the first to address potential failure through
approximate second-order points, which as highlighted in \exaref{approx},
can occur even when spurious local minima do not exist. Through a
\emph{sharp} characterization of approximate second-order points in
the noisy setting, we discover that near-second-order points of the
slightly overparameterized problem with $r=\Omega(r^{\star})$ achieve
the same minimax-optimal error rates as classical convex approaches.
The new ability to accommodate approximate second-order points also
provides us with polynomial-time complexity guarantees.

Finally, we mention that a slight \emph{additive} overparameterization
$r-r^{\star}=\Omega(1)$ was recently shown to induce benign landscape
in the phase and orthogonal group synchronization problems~\cite{bandeira2016low,ling2023solving,mcrae2024benign}.
This is quite a different setting to ours, but the additive guarantee
is also more favorable than the multiplicative one $r/r^{\star}=\Omega(1)$
that we show to be necessary and sufficient in our setting. It would
therefore be interesting to see whether our sharp multiplicative guarantee
can be improved to an additive one, e.g. by imposing additional structure
of the measurements, or by adopting further regularization. 

\subsection{Organization}

Our main contribution is the new proof technique for establishing
sufficient conditions in the overparameterized regime. In \secref{duality},
we begin by outlining the strong duality between the existing proof
techniques, as it is the central idea behind our proof technique.
In \secref{proof}, we state the sharp escape directions corresponding
to \thmref{zhang}, and use these to prove our new sufficient conditions.
In \secref{counter}, we state the counterexamples used to prove our
new necessary conditions. Finally, in \secref{conclusion}, we offer
some concluding remarks on extending our analysis to non-RIP cases,
such as the matrix completion problem. 

\subsection*{Notation}

Our notation is  consistent with MATLAB syntax. We use the comma
for horizontal concatenation $[a,b]=\left[\begin{smallmatrix}a & b\end{smallmatrix}\right]$,
and the semicolon for vertical concatenation $[a;b]=\left[\begin{smallmatrix}a\\
b
\end{smallmatrix}\right]$. For $w\in\R^{n}$, denote $\diag(w)$ as the $n\times n$ diagonal
matrix with $w$ along its diagonal. For $W\in\R^{n\times n}$, denote
$\diag(W)$ as the length-$n$ vector constructed by indexing the
diagonal elements of $W$. 

\section{\protect\label{sec:duality}Strong duality between counterexamples
and escape directions}

In the RIP setting, \thmref{zhang} was the first to uncover an unambiguous
improvement with overparameterization. This was achieved by fully
sharpening the conservatism of prior global guarantees using a novel
proof technique, based on the \emph{inexistence of a counterexample},
which was originally introduced by \cite{zhang2019sharp}. To prove
that all functions of the following class have no spurious local minima
\[
\cF_{Z}(\delta,\epsilon_{0})=\left\{ f:\R^{n\times r}\to\R:\begin{array}{c}
f(X)=\frac{1}{2}\|\AA(XX^{T})-b\|^{2},\\
\|\AA^{T}(\AA(ZZ^{T})-b)\|_{\op}\le\epsilon_{0},\\
\AA\in\ripd
\end{array}\right\} 
\]
the basic idea is to demonstrate that a counterexample to refute this
claim does not exist. The motivating insight behind this approach,
first noted in~\cite{zhang2018much}, is that the search for a counterexample
with a \emph{fixed} spurious point $X\in\R^{n\times r}$ and \emph{fixed}
ground truth $Z\in\R^{n\times r^{\star}}$, as follows 
\begin{alignat}{2}
\text{find }\quad & f\in\cF_{Z}(\delta,\epsilon_{0}) &  & \tag{P}\label{eq:Pnoisy}\\
\text{ such that }\quad & \inner{\nabla f(X)}D\ge-\epsilon_{1}\|DX^{T}\| &  & \text{for all }D\in\R^{n\times r},\nonumber \\
 & \inner{\nabla f(X)[D]}D\ge-\epsilon_{2}\|DX^{T}\|^{2} & \quad & \text{for all }D\in\R^{n\times r},\nonumber 
\end{alignat}
can be formulated and solved as a convex semidefinite program (SDP).
In the special case of $\epsilon_{0}=\epsilon_{1}=\epsilon_{2}=0$,
corresponding to exact recovery from noiseless measurements with exact
second-order points, Zhang \cite{zhang2022improved} exhaustively
solved all instances of (\ref{eq:Pnoisy}) in closed form, across
all possible spurious points $X$ and ground truths $Z$ where $XX^{T}\ne ZZ^{T}$.
This way, he found that counterexamples exist if and only if the RIP
constant satisfies $\delta\ge(1+\sqrt{r^{\star}/r})^{-1}$. The inexistence
of a counterexample $f\in\cF_{Z}(\delta,\epsilon_{0})$ with $\delta<(1+\sqrt{r^{\star}/r})^{-1}$
thus proves that spurious local minima do not exist, so exact recovery
is guaranteed. 

Unfortunately, the need to exhaust all possible counterexamples also
makes the proof technique very fragile and difficult to generalize.
While it is easy to look for a counterexample, by solving (\ref{eq:Pnoisy})
numerically across many choices of $X$ and $Z$, failing to find
one does not prove that they do not exist. To rigorously rule out
all possible counterexamples, one would need to solve all instances
of (\ref{eq:Pnoisy}) in closed form. But Zhang's closed-form solution
for $\epsilon_{0}=\epsilon_{1}=\epsilon_{2}=0$ is already very complicated,
and it is unclear how it can be further generalized to noisy measurements
$\epsilon_{0}>0$ and approximate second-order points $\epsilon_{1},\epsilon_{2}>0$,
or whether a closed-form solution even exists. Indeed, this inability
to solve more general instances of (\ref{eq:Pnoisy}) in closed form
was cited as a major obstacle in prior applications of the proof technique~\cite{zhang2020many,ma2023geometric}. 

In this paper, we take the dual of the above approach; we state Lagrange
multipliers $D_{0},D_{1},\dots,D_{nr}$ for each pair of spurious
point $X$ and ground truth $Z$ that solve the following separating
hyperplane problem
\begin{equation}
\text{find }\quad D_{0},D_{1},\dots,D_{nr}\quad\text{ such that }\quad\psi(D_{0},D_{1},\dots,D_{nr})<0,\tag{D}\label{eq:Dnoisy}
\end{equation}
in which the dual function $\psi$ is defined
\[
\psi(D_{0},D_{1},\dots,D_{nr})\eqdef\sup_{f\in\cF_{Z}(\delta,\epsilon_{0})}\left\{ \begin{array}{c}
(\inner{\nabla f(X)}{D_{0}}+\epsilon_{1}\|D_{0}X^{T}\|)\\
+\frac{1}{2}\sum_{i=1}^{nr}(\inner{\nabla^{2}f(X)[D_{i}]}{D_{i}}+\epsilon_{2}\|D_{i}X^{T}\|^{2})
\end{array}\right\} .
\]
It follows from basic Lagrangian duality that feasibility in (\ref{eq:Dnoisy})
certifies infeasibility in (\ref{eq:Pnoisy}). In fact, this is
the exact approach based on the \emph{existence of an escape direction}
used to prove the majority of existing benign landscape guarantees~\cite{bhojanapalli2016global,ge2017nospurious,park2017non,zhu2018global,li2019non}.
Typically, the dual solutions $D_{0},D_{1},\dots,D_{nr}$ are constructed
using the displacement vector~\cite{zheng2015convergent,tu2016low}
\[
\Delta=Z_{X}-X\qquad\text{ where }Z_{X}=\arg\min_{Z_{X}Z_{X}^{T}=ZZ^{T}}\|Z_{X}-X\|.
\]
For example, in their elegantly simple proof, \cite{ge2017nospurious}
used $D_{0}=4\Delta,$ $D_{1}=\Delta,$ and $D_{2}=\cdots=D_{nr}=0$
to prove for the special case of $\epsilon_{0}=\epsilon_{1}=\epsilon_{2}=0$
that
\begin{align*}
\sup_{f\in\cF_{Z}(\delta)}\left\{ \inner{\nabla f(X)}{4\Delta}+\inner{\nabla^{2}f(X)[\Delta]}{\Delta}\right\} \le & -(1-5\delta)\|XX^{T}-ZZ^{T}\|^{2}.
\end{align*}
This immediately implies for $f\in\cF_{Z}(\delta,\epsilon_{0})$ with
$\delta<1/5$ that spurious local minima do not exist, so exact recovery
is guaranteed. 

The main feature of the dual approach is that it is easy to generalize;
once a dual solution $D_{0},D_{1},\dots,D_{nr}$ has been identified
for the one setting, it can be quickly reused for another setting,
such as the asymmetric parameterization~\cite{park2017non}, and
even non-RIP recovery like the matrix completion problem~\cite{ge2017nospurious}.
 Unfortunately, none of the existing guarantees proved using the
dual approach have been sharp enough to exhibit a dependence on the
overparameterization ratio $r/r^{\star}$. A critical concern is that
the dual approach might be fundamentally too conservative. For a \emph{nonconvex}
feasibility problem like (\ref{eq:Pnoisy}), it is possible---and
indeed, common---for the primal to be infeasible without there existing
a corresponding dual solution to certify this fact.  Surprisingly,
we prove in this paper that strong duality holds between (\ref{eq:Pnoisy})
and (\ref{eq:Dnoisy}); the inexistence of a counterexample in (\ref{eq:Pnoisy})
can \emph{always} be certified by the existence of  escape directions
$D_{0},D_{1},\dots,D_{nr}$ in (\ref{eq:Dnoisy}). 
\begin{theorem}[Strong duality]
\label{thm:strongdual}Fix $X\in\R^{n\times r}$ and $Z\in\R^{n\times r^{\star}}$
such that $r\ge r^{\star}$ and $XX^{T}\ne ZZ^{T}$. For RIP constant
$0\le\delta<1$ and accuracy parameters $\epsilon_{0},\epsilon_{1},\epsilon_{2}\ge0$,
exactly one of the following statements is true:
\begin{itemize}
\item (Counterexample) There exists $f\in\cF_{Z}(\delta,\epsilon_{0})$
for which $X$ is an approximate spurious second-order point: 
\[
\inner{\nabla f(X)}D\ge-\epsilon_{1}\|DX^{T}\|,\quad\inner{\nabla^{2}f(X)[D]}D\ge-\epsilon_{2}\|DX^{T}\|^{2}\quad\text{for all }D\in\R^{n\times r}.
\]
\item (Escape directions) There exists $D_{0},D_{1},\dots,D_{nr}$ that
together guarantee a local decrement at $X$ across all arbitrary
$f\in\cF_{Z}(\delta,\epsilon_{0})$:
\[
\sup_{f\in\cF_{Z}(\delta,\epsilon_{0})}\left\{ \inner{\nabla f(X)}{D_{0}}+\epsilon_{1}\|D_{0}X^{T}\|+\sum_{i=1}^{nr}\left(\inner{\nabla^{2}f(X)[D_{i}]}{D_{i}}+\epsilon_{2}\|D_{i}X^{T}\|^{2}\right)\right\} <0.
\]
\end{itemize}
\end{theorem}

Therefore, the dual approach comes with no loss of sharpness compared
to the primal approach. In particular, \thmref{strongdual} guarantees
the existence of sharp escape directions that would improve existing
escape-direction proofs to the same sharpness as \thmref{zhang}.
In \secref{proof}, we explicitly identify these sharp escape directions,
and use these to prove our main results \thmref{sym} and \thmref{asym}. 

In the rest of this section, we give a proof of \thmref{strongdual}.
Our proof works by appealing to an equivalent statement posed over
the following family of quadratic functions
\[
\cH_{Z}(\delta,\epsilon_{0})=\left\{ h:\R^{n\times r}\to\R:\begin{array}{c}
h(U)=\inner G{UU^{T}}+\half\|UU^{T}-ZZ^{T}\|_{\HH}^{2},\\
\|G\|_{\op}\le\epsilon_{0},\quad(1-\delta)\Id\preceq\HH\preceq(1+\delta)\Id.
\end{array}\right\} 
\]
Critically, the set $\cH_{Z}(\delta,\epsilon_{0})$ is convex: if
$h_{1},h_{2}\in\cH_{Z}(\delta,\epsilon_{0})$, then $th_{1}+(1-t)h_{2}\in\cH_{Z}(\delta,\epsilon_{0})$
for all $0\le t\le1$. Therefore, we can generally expect strong duality
to hold for existential statements posed over $\cH_{Z}(\delta,\epsilon_{0})$.
\begin{lemma}
\label{lem:cvxdual}Under the same setting as \thmref{strongdual},
exactly one of the following is true:
\begin{itemize}
\item (Counterexample) There exists $h\in\cH_{Z}(\delta,\epsilon_{0})$
such that
\[
\inner{\nabla h(X)}D\ge-\epsilon_{1}\|DX^{T}\|,\quad\inner{\nabla^{2}f(X)[D]}D\ge-\epsilon_{2}\|DX^{T}\|^{2}\quad\text{for all }D\in\R^{n\times r}.
\]
\item (Escape directions) There exists $D_{0},D_{1},\dots,D_{nr}$ such
that
\[
\sup_{h\in\cH_{Z}(\delta,\epsilon_{0})}\left\{ \inner{\nabla h(X)}{D_{0}}+\epsilon_{1}\|D_{0}X^{T}\|+\sum_{i=1}^{nr}\left(\inner{\nabla^{2}h(X)[D_{i}]}{D_{i}}+\epsilon_{2}\|D_{i}X^{T}\|^{2}\right)\right\} <0.
\]
\end{itemize}
\end{lemma}

\begin{proof}
The amounts to verifying that the dual problem satisfies Slater's
condition. 
\end{proof}
We begin by attempting to construct a counterexample $f\in\cF_{Z}(\delta,\epsilon_{0})$
with ground truth $Z$ and spurious second-order point $X$, by constructing
a counterexample $h\in\cH_{\hZ}(\delta,\epsilon_{0})$ with ground
truth $\hZ=P^{T}Z$ and spurious second-order point $\hX=P^{T}X$,
where $P=\orth([X,Z])$. If this suceeds, then we can evoke the following
lemma to fulfill the first clause of \thmref{strongdual}.
\begin{lemma}
Given $X\in\R^{n\times r}$ and $Z\in\R^{n\times r^{\star}}$, define
$\hX=P^{T}X\in\R^{d\times r}$ and $\hat{Z}=P^{T}Z\in\R^{d\times r^{\star}}$
where $P=\orth([X,Z])\in\R^{n\times d}$. For every $h\in\cH_{\hZ}(\delta,\epsilon_{0})$,
there exists a choice of $f\in\cF_{Z}(\delta,\epsilon_{0})$ such
that 
\[
\inner{\nabla f(X)}D=\inner{\nabla h(\hX)}{P^{T}D},\;\inner{\nabla^{2}f(X)[D]}D=\inner{\nabla h(\hX)[P^{T}D]}{P^{T}D}+2\|P_{\perp}^{T}D\hX^{T}\|^{2}
\]
where $P_{\perp}\in\R^{n\times(n-d)}$ denotes the orthogonal complement
of $P$.
\end{lemma}

\begin{proof}
Given $h(U)=\inner G{UU^{T}}+\frac{1}{2}\|UU^{T}-ZZ^{T}\|_{\HH}^{2}$,
we factor $\HH=\UU^{T}\UU$ into its upper-triangular Cholesky factor
$\UU:\S^{d}\to\R^{\frac{1}{2}d(d+1)}$, and define
\begin{align*}
\AA(M) & =(\half\UU(M_{11}+M_{11}^{T}),\half\vect(M_{11}-M_{11}^{T}),\vect(M_{21}),\vect(M_{12}),\vect(M_{22}))
\end{align*}
where $M_{11}=P^{T}MP,$ $M_{12}=P^{T}MP_{\perp},$ $M_{21}=P_{\perp}^{T}MP,$
and $M_{22}=P_{\perp}^{T}MP_{\perp}$. We can verify that $(1-\delta)\Id\preceq\AA^{T}\AA\preceq(1+\delta)\Id$
holds, so $\mathrm{Im}(\AA^{T})=\R^{n\times n}$, and there must exist
$b$ such that $\AA^{T}[\AA(ZZ^{T})-b]=PGP^{T}$. Then, $f(U)=\frac{1}{2}\|\AA(UU^{T})-b\|^{2}$
satisfies $f\in\cF_{Z}(\delta,\epsilon_{0})$, because $\AA\in\rip(\delta,n)\subseteq\ripd$
and $\|\AA^{T}[\AA(ZZ^{T})-b]\|_{\op}=\|G\|_{\op}\le\epsilon_{0}$.
Finally, for this choice of $f$, we can mechanically verify that
$\inner{\nabla f(X)}D=\inner{\nabla h(\hX)}{P^{T}D}$ and $\inner{\nabla^{2}f(X)[D]}D=\inner{\nabla h(\hX)[P^{T}D]}{P^{T}D}+2\|P_{\perp}^{T}D\hX^{T}\|^{2}$.
\end{proof}
Conversely, if a counterexample $h\in\cH_{\hZ}(\delta,\epsilon_{0})$
with ground truth $\hZ=P^{T}Z$ and spurious second-order point $\hX=P^{T}X$
does not exist, then strong duality in \lemref{cvxdual} guarantees
the existence of escape directions $\hD_{0},\hD_{1},\dots,\hD_{nr}$
that yield a local decrement across all $h\in\cH_{\hZ}(\delta,\epsilon_{0})$
at $\hX$. Then, we can evoke the following lemma to show that $D_{0},D_{1},\dots,D_{nr}$
with $D_{i}=P\hD_{i}$ are escape directions that escape all $f\in\cF_{Z}(\delta,\epsilon_{0})$
at $X$, hence fulfilling the second clause of \thmref{strongdual}. 
\begin{lemma}
Given $X\in\R^{n\times r}$ and $Z\in\R^{n\times r^{\star}}$, define
$\hX=P^{T}X\in\R^{d\times r}$ and $\hat{Z}=P^{T}Z\in\R^{d\times r^{\star}}$
where $P=\orth([X,Z])\in\R^{n\times d}$. For every $f\in\cF_{Z}(\delta,\epsilon_{0})$,
there exists a choice of $h\in\cH_{\hZ}(\delta,\epsilon_{0})$ such
that 
\[
\inner{\nabla h(\hX)}{\hat{D}}=\inner{\nabla f(X)}{P\hat{D}},\quad\inner{\nabla^{2}h(\hX)[\hat{D}]}{\hat{D}}=\inner{\nabla f(X)[P\hat{D}]}{P\hat{D}}.
\]
\end{lemma}

\begin{proof}
Given $f(U)=\frac{1}{2}\|\AA(UU^{T})-b\|^{2}$, we define $\HH$ such
that $\|M\|_{\HH}=\|\AA(PMP^{T})\|^{2}$ for all $M\in\S^{d}$, and
$G=P^{T}[\AA^{T}(\AA(ZZ^{T})-b)]P$. Then, $h(\hU)=\inner G{\hU\hU^{T}}+\frac{1}{2}\|\hU\hU^{T}-\hZ\hZ^{T}\|^{2}$
satisfies $h\in\cH_{\hZ}(\delta,\epsilon_{0})$, because $\AA\in\ripd$
and $\rank(PMP^{T})\le d\le r+r^{\star}$ together imply $\left|\|M\|_{\HH}^{2}/\|M\|^{2}-1\right|\le\delta$
for all $M\in\R^{d}$. Finally, for this choice of $h$, we can mechanically
verify that $\inner{\nabla h(\hX)}{\hat{D}}=\inner{\nabla f(X)}{P\hat{D}}$
and $\inner{\nabla^{2}h(\hX)[\hat{D}]}{\hat{D}}=\inner{\nabla f(X)[P\hat{D}]}{P\hat{D}}$.
\end{proof}

\section{\protect\label{sec:proof}Proof of sufficiency by sharp escape directions}

The main technical contribution of this paper is to explicitly identify
the sharp escape directions stated below. 
\begin{definition}[Sharp escape directions]
\label{def:directions}For $X\in\R^{n\times r}$ and $Z\in\R^{n\times r^{\star}}$
such that $r\ge r^{\star}$. Define $D_{0},D_{1},\dots,D_{r^{\star}}\in\R^{n\times r}$
as follows
\[
D_{0}=(I-\half XX^{\dagger})Z(X^{\dagger}Z)^{T}-\half X,\qquad D_{i}=Z_{\perp}e_{i}u^{T}\quad\text{ for }i\in\{1,2,\dots,r^{\star}\},
\]
where $Z_{\perp}=(I-XX^{\dagger})Z$ and $u\in\R^{r},\|u\|=1$ is
the eigenvector satisfying $u^{T}X^{T}Xu=\lambda_{\min}(X^{T}X)$.
\end{definition}

These are chosen explicitly to decompose the error vector $E=XX^{T}-ZZ^{T}$
into orthogonal components $E_{\T}$ and $E_{\N}$ that lie respectively
in the tangent and normal spaces of the Riemannian manifold of rank-$r$
positive semidefinite matrices.
\begin{fact}[Tangent-normal decomposition]
\label{fact:decomp}Under the same setting as \defref{directions},
the error vector $E=XX^{T}-ZZ^{T}$ decomposes into orthogonal components
$E=E_{\T}+E_{\N}$ 
\begin{alignat*}{2}
E_{\T} & =-(D_{0}X^{T}+XD_{0}^{T}) & \qquad & =\Pi_{X}E+E\Pi_{X}-\Pi_{X}E\Pi_{X},\\
E_{\N} & =-\sum_{i=1}^{r^{\star}}D_{i}D_{i}^{T} &  & =(I-\Pi_{X})E(I-\Pi_{X})=-Z_{\perp}Z_{\perp}^{T},
\end{alignat*}
where $\Pi_{X}=XX^{\dagger}$ is the projector onto the column space
of $X$.
\end{fact}

Of all choices of $D_{0},D_{1},\dots,D_{r^{\star}}$ that decompose
the error into a tangent and orthogonal component, our choice is made
to minimize the combined local norm $\sum_{i=0}^{r^{\star}}\|D_{i}X^{T}\|^{2}$. 
\begin{fact}[Minimum local norm]
\label{fact:nrm}The search directions $D_{0},D_{1},\dots,D_{r^{\star}}$
satisfy
\begin{align*}
\|D_{0}X^{T}\| & =\frac{1}{2}\|D_{0}X^{T}+XD_{0}^{T}\|,\\
\sum_{i=1}^{r^{\star}}\|D_{i}X^{T}\|^{2} & =\frac{1}{2}\sum_{i=1}^{r^{\star}}\|D_{i}X^{T}+XD_{i}^{T}\|^{2}=\lambda_{\min}(X^{T}X)\|Z_{\perp}\|^{2}.
\end{align*}
\end{fact}

Before we proceed with the proof of our main results, we first establish
two small technical lemmas regarding this tangent-normal decomposition.
Below and henceforth, we use the pseudo-norm notation $\|M\|_{\HH}=\sqrt{\inner{\HH(M)}M}$
for some self-adjoint linear operator $\HH$, without necessarily
requiring $\HH$ to be positive definite.
\begin{lemma}
Under the same setting as \factref{decomp}, if $\HH:\S^{n}\to\S^{n}$
satisfies $\left|\|M\|_{\HH}^{2}/\|M\|^{2}-1\right|\le\delta$ for
all $M\in\S^{n}$ such that $\rank(M)\le r+r^{\star}$, then it also
satisfies
\[
\inner{\HH(E)}F\ge\inner EF-\delta\|E\|\|F\|\quad\text{for all }F=t_{1}E_{\T}+t_{2}E_{\N},\;t_{1},t_{2}\in\R.
\]
\end{lemma}

\begin{proof}
Without loss of generality, let $\|E\|=\|F\|=1$. Define $\cM_{P}=\{P\hat{M}P^{T}:\hat{M}\in\S^{d}\}$
where $P=\orth([X,Z])\in\R^{n\times d}$. Observe that $E\in\cM_{P}$
and $E_{\N}\in\cM_{P}$, and $E_{\T}=E-E_{\N}\in\cM_{P}$ by the convexity
of $\cM_{P}$. Therefore, $E+F\in\cM_{P}$ and $E-F\in\cM_{P}$ because
$F=t_{1}E_{\T}+t_{2}E_{\N}$. It now follows from the fact that $\rank(M)\le r+r^{\star}$
for every $M\in\cM_{P}$ that 
\begin{align*}
4\inner{\HH(E)}F & =\|E+F\|_{\HH}^{2}-\|E-F\|_{\HH}^{2}\ge(1-\delta)\|E+F\|^{2}-(1+\delta)\|E-F\|^{2}\\
 & =4\inner EF-2\delta(\|E\|^{2}+\|F\|^{2})=4\inner EF-4\delta.
\end{align*}
\end{proof}
\begin{lemma}
\label{lem:slacklb}Under the same setting as \factref{decomp}, if
the matrix $N\succeq0$ satisfies $NZ=0$, then 
\[
\inner NE\ge0,\qquad\inner N{t_{1}E_{\T}+t_{2}E_{\N}}\ge0\text{ for all }t_{1}\ge t_{2}\ge0.
\]
\end{lemma}

\begin{proof}
We have $\inner NE=\inner N{XX^{T}-ZZ^{T}}=\inner N{XX^{T}}\ge0,$
so
\begin{gather*}
\inner N{t_{1}E_{\T}+t_{2}E_{\N}}=\inner N{t_{1}E-(t_{1}-t_{2})E_{\N}}=\inner N{t_{1}XX^{T}+(t_{1}-t_{2})Z_{\perp}Z_{\perp}^{T}}\ge0.
\end{gather*}
\end{proof}
We are now ready to prove the sufficient conditions in our main results.
In \secref{sym}, we provide a complete, detailed proof of the symmetric
case (\thmref{sym} and \corref{sym}). Then, in \secref{asym}, we
describe how this proof readily generalizes to the asymmetric case
(\thmref{asym} and \corref{nonsym}). 

\subsection{\protect\label{sec:sym}Symmetric parameterization $XX^{T}$}

Our proof of \thmref{sym} follows by plugging a rescaled version
of $D_{0},D_{1},\dots,D_{r^{\star}}$ into the dual problem (\ref{eq:Dnoisy}).
Concretely, for every spurious point $X$ and ground truth $Z$ with
$XX^{T}\ne ZZ^{T}$, we identify a rescaling $t_{1}\ge0,t_{2}\ge0$
so that $\psi(t_{1},t_{2})<0$ holds where
\begin{equation}
\psi(t_{1},t_{2})\eqdef\sup_{f\in\cF_{Z}(\delta,\epsilon_{0})}\left\{ \begin{array}{c}
t_{1}(\inner{\nabla f(X)}{D_{0}}+\epsilon_{1}\|D_{0}X^{T}\|)\\
+\half t_{2}\sum_{i=1}^{nr}(\inner{\nabla^{2}f(X)[D_{i}]}{D_{i}}+\epsilon_{2}\|D_{i}X^{T}\|^{2})
\end{array}\right\} .\label{eq:dual1}
\end{equation}
Concrete expressions for the directional derivatives of $f\in\cF_{Z}(\delta,\epsilon_{0})$
are stated below; we recall that $f(X)=\frac{1}{2}\|\AA(XX^{T})-b\|^{2}$
with noisy measurements $b$ satisfying $\|\AA^{T}(b-\AA(ZZ^{T}))\|_{\op}\le\epsilon_{0}$.
 
\begin{fact}
\label{fact:deriv1}For $f\in\cF_{Z}(\delta,\epsilon_{0})$, write
$S=\AA^{T}(\AA(ZZ^{T})-b)$ and $\HH(M)\eqdef\AA^{T}\AA(M)$. Then,
the directional derivatives of $f$ are written
\begin{align*}
\inner{\nabla f(X)}D & =\inner{\HH(E)+S}{DX^{T}+XD^{T}},\\
\half\inner{\nabla^{2}f(X)[D]}D & =\inner{\HH(E)+S}{DD^{T}}+\frac{1}{2}\|DX^{T}+XD^{T}\|_{\HH}^{2}.
\end{align*}
Moreover, $\|S\|_{\op}\le\epsilon_{0}$ and $\left|\|M\|_{\HH}^{2}/\|M\|^{2}-1\right|\le\delta$
for all $M\in\S^{n}$ such that $\rank(M)\le r+r^{\star}$.
\end{fact}

We first consider two trivial cases, which correspond to cases where
the error vector $E$ lies entirely in the tangent or the normal spaces.
\begin{lemma}
\label{lem:trivial1}Fix $X\in\R^{n\times r}$ and $Z\in\R^{n\times r^{\star}}$
such that $r\ge r^{\star}$ and $E=XX^{T}-ZZ^{T}\ne0$. Define the
dual function $\psi(t_{1},t_{2})$ as in (\ref{eq:dual1}). If $X\ne0$
and $XX^{\dagger}Z=Z$, then 
\[
\psi(1,0)\le(\sqrt{r+r^{\star}}\epsilon_{0}+\half\epsilon_{1})\|E\|-(1-\delta)\|E\|^{2}.
\]
If instead $X=0$, and hence $XX^{\dagger}Z=0$, then 
\[
\psi(0,1)\le\sqrt{r+r^{\star}}\epsilon_{0}\|E\|-(1-\delta)\|E\|^{2}.
\]
\end{lemma}

\begin{proof}
Decompose the error $E=E_{\T}+E_{\N}$ as in \factref{decomp}. If
$X\ne0$ and $(I-XX^{\dagger})Z=Z_{\perp}=0$, then $E_{\N}=-Z_{\perp}Z_{\perp}^{T}=0$
and $E_{\T}=-(D_{0}X^{T}+XD_{0}^{T})=E$. Hence, 
\begin{align*}
\sup_{f\in\cF_{Z}(\delta,\epsilon_{0})}\inner{\nabla f(X)}{D_{0}}+\epsilon_{1}\|D_{0}X^{T}\|= & \sup_{f\in\cF_{Z}(\delta,\epsilon_{0})}\inner{\HH(E)+S}{-E_{\T}}+\frac{\epsilon_{1}}{2}\|E_{\T}\|\\
= & \sup_{f\in\cF_{Z}(\delta,\epsilon_{0})}-\|E\|_{\HH}^{2}-\inner SE+\frac{\epsilon_{1}}{2}\|E\|\\
\le & -(1-\delta)\|E\|^{2}+\epsilon_{0}\|E\|_{\nuc}+\frac{\epsilon_{1}}{2}\|E\|.
\end{align*}
Similarly, if $X=0$, then $E_{T}=-(D_{0}X^{T}+XD_{0}^{T})=0$ and
$E_{\N}=-\sum_{i>0}D_{i}D_{i}^{T}=E$. It follows from $\|D_{i}X^{T}\|=0$
that
\begin{align*}
\sup_{f\in\cF_{Z}(\delta,\epsilon_{0})}\frac{1}{2}\sum_{i=1}^{r^{\star}}\left(\inner{\nabla^{2}f(X)[D_{i}]}{D_{i}}+\epsilon_{2}\|D_{i}X^{T}\|^{2}\right) & =\sup_{f\in\cF_{Z}(\delta,\epsilon_{0})}\inner{\HH(E)+S}{-E_{N}}\\
 & =\sup_{f\in\cF_{Z}(\delta,\epsilon_{0})}-\|E\|_{\HH}^{2}-\inner SE\\
 & \le-(1-\delta)\|E\|^{2}+\epsilon_{0}\|E\|_{\nuc}.
\end{align*}
Finally, in both cases, we bound $\|E\|_{\nuc}\le\sqrt{\rank(E)}\|E\|$. 
\end{proof}
Next, we consider the nontrivial cases where the error vector $E$
has components in both the tangent and normal spaces.

\begin{lemma}
\label{lem:nontriv1}Fix $X\in\R^{n\times r}$ and $Z\in\R^{n\times r^{\star}}$
such that $r\ge r^{\star},$ $E=XX^{T}-ZZ^{T}\ne0,$ $X\ne0,$ and
$Z_{\perp}=(I-XX^{\dagger})\ne0$. Define the parameters $\alpha,\beta$
as follows 
\[
\alpha=\frac{\|Z_{\perp}Z_{\perp}^{T}\|}{\|E\|},\qquad\beta=\frac{\lambda_{\min}(X^{T}X)}{\|E\|}\frac{\tr(Z_{\perp}Z_{\perp}^{T})}{\|Z_{\perp}Z_{\perp}^{T}\|}.
\]
Then, the dual function $\psi(t_{1},t_{2})$ in (\ref{eq:dual1})
is upper-bounded
\begin{align*}
\psi(\sqrt{\frac{1-\tau^{2}}{1-\alpha^{2}}},\frac{\tau}{\alpha}) & \le(\epsilon_{0}\sqrt{r+r^{\star}}+\half\epsilon_{1})\|E\|\\
 & \qquad-\left[\sqrt{1-\alpha^{2}}\sqrt{1-\tau^{2}}+[\alpha-(1+\delta+\half\epsilon_{2})\beta]\tau-\delta\right]\|E\|^{2}.
\end{align*}
\end{lemma}

\begin{proof}
Write $t_{1}=\sqrt{\frac{1-\tau^{2}}{1-\alpha^{2}}}$ and $t_{2}=\frac{\tau}{\alpha}$
and $F=t_{1}E_{\T}+t_{2}E_{N}$. We repeat the definitions and manipulations
in the proof of \lemref{trivial1} to obtain the following for $\epsilon_{1}=\epsilon_{2}=0$:
\begin{align*}
\psi_{0}(t_{1},t_{2})\eqdef & \sup_{f\in\cF_{Z}(\delta,\epsilon_{0})}\quad t_{1}\inner{\nabla f(X)}{D_{0}}+\frac{t_{2}}{2}\sum_{i=1}^{r^{\star}}\inner{\nabla^{2}f(X)[D_{i}]}{D_{i}}\\
= & \sup_{f\in\cF_{Z}(\delta,\epsilon_{0})}\quad-\inner{\HH(E)+S}F+\frac{t_{2}}{2}\sum_{i=1}^{r^{\star}}\|D_{i}X^{T}+XD_{i}^{T}\|_{\HH}^{2},\\
\le & -(\inner EF-\delta\|E\|\|F\|)+\epsilon_{0}\|F\|_{\nuc}+t_{2}(1+\delta)\lambda_{\min}(X^{T}X)\|Z_{\perp}\|^{2}.
\end{align*}
Then, the general case of $\epsilon_{1}\ge0,\epsilon_{2}\ge0$ is
simply
\begin{align*}
\psi(t_{1},t_{2}) & =\psi_{0}(t_{1},t_{2})+t_{1}\epsilon_{1}\|D_{0}X^{T}\|+\frac{t_{2}\epsilon_{2}}{2}\sum_{i=1}^{r^{\star}}\|D_{i}X^{T}\|^{2}\\
 & =\psi_{0}(t_{1},t_{2})+\frac{t_{1}\epsilon_{1}}{2}\|E_{\T}\|+\frac{t_{2}\epsilon_{2}}{2}\lambda_{\min}(X^{T}X)\|Z_{\perp}\|^{2}.
\end{align*}
Finally, we verify that
\begin{gather*}
\inner EF=\inner E{t_{1}E_{\T}+t_{2}E_{N}}=t_{1}\|E_{\T}\|^{2}+t_{2}\|E_{\N}\|^{2}=(\sqrt{1-\alpha^{2}}\sqrt{1-\tau^{2}}+\alpha\tau)\|E\|^{2},\\
\|F\|^{2}=\|t_{1}E_{\T}+t_{2}E_{N}\|^{2}=t_{1}^{2}\|E_{\T}\|^{2}+t_{2}^{2}\|E_{\N}\|^{2}=\|E\|^{2},\\
\|F\|_{\nuc}\le\sqrt{\rank(F)}\|F\|\le\sqrt{r+r^{\star}}\|F\|,\\
t_{2}\lambda_{\min}(X^{T}X)\|Z_{\perp}\|^{2}=t_{2}\alpha\beta\|E\|^{2}=\tau\beta\|E\|^{2}.
\end{gather*}

\end{proof}

The optimal choice of $\tau$ for the nontrivial case above is determined
through the following lemma, which at its heart is a purely linear
algebraic result. 
\begin{lemma}
\label{lem:techical}Under the same setting as \lemref{nontriv1},
we have

\[
\max_{0\le\tau\le\alpha}\left\{ \frac{\sqrt{1-\alpha^{2}}\sqrt{1-\tau^{2}}+\alpha\tau-(1+\frac{1}{2}\epsilon_{2})\beta\tau}{1+(1+\frac{1}{2}\epsilon_{2})\beta\tau}\right\} \ge\frac{1}{1+(1+\frac{1}{2}\epsilon_{2})\sqrt{r^{\star}/r}}.
\]
\end{lemma}

\begin{proof}
\global\long\def\tbeta{\tilde{\beta}}%
We begin by citing two important lemmas from Zhang's~\cite{zhang2022improved}
proof of \thmref{zhang}. First, the following problem has closed-form
solution~\cite[Lemma~5.4]{zhang2022improved}
\[
\delta(\alpha,\beta)\eqdef\max_{0\le\tau\le\alpha}\left\{ \frac{\sqrt{1-\alpha^{2}}\sqrt{1-\tau^{2}}+\alpha\tau-\beta\tau}{1+\beta\tau}\right\} =\begin{cases}
\sqrt{1-\alpha^{2}} & \beta\ge\frac{\alpha}{1+\sqrt{1-\alpha^{2}}},\\
\frac{1-2\alpha\beta+\beta^{2}}{1-\beta^{2}} & \beta\le\frac{\alpha}{1+\sqrt{1-\alpha^{2}}}.
\end{cases}
\]
Second, the two parameters $\alpha,\beta$ as defined in \lemref{nontriv1}
will always satisfy $\alpha^{2}+(r/r^{\star})\beta^{2}\le1-[(\beta-\alpha)_{+}]^{2}$~\cite[Lemma~3.8]{zhang2022improved}.
This motivates us to use the following as a lower-bound
\[
\delta^{\star}=\min_{\alpha\ge0,\beta\ge0}\{\delta(\alpha,(1+\half\epsilon_{2})\beta):\alpha^{2}+\frac{r}{r^{\star}}\beta^{2}\le1-[(\beta-\alpha)_{+}]^{2}\}.
\]
Write $\tbeta=(1+\half\epsilon_{2})\beta$. First, we consider the
region $\tbeta\ge{\displaystyle \frac{\alpha}{1+\sqrt{1-\alpha^{2}}}}$.
We solve the following for every fixed $\tilde{\beta}$, and find
that the minimum is attained at the boundary 
\[
\min_{\alpha\ge0}\left\{ \sqrt{1-\alpha^{2}}:\tbeta\ge{\displaystyle \frac{\alpha}{1+\sqrt{1-\alpha^{2}}}}\right\} =\frac{1-2\alpha\tbeta+\tbeta^{2}}{1-\tbeta^{2}}.
\]
Indeed, substituting $\tbeta=\frac{\alpha}{1+\sqrt{1-\alpha^{2}}}=\frac{1+\sqrt{1-\alpha^{2}}}{\alpha}$
into the following yields
\[
\left(\frac{1-2\alpha\tbeta+\tbeta^{2}}{1-\tbeta^{2}}\right)\left(\frac{\alpha/\tbeta}{\alpha/\tbeta}\right)=\frac{\alpha/\tbeta-2\alpha^{2}+\alpha\beta}{\alpha/\beta-\alpha\beta}=\frac{2-2\alpha^{2}}{2\sqrt{1-\alpha^{2}}}=\sqrt{1-\alpha^{2}}.
\]
Next, we consider the region $\tilde{\beta}\le{\displaystyle \frac{\alpha}{1+\sqrt{1-\alpha^{2}}}}$.
Note that $\alpha\ge\tilde{\beta}\ge\beta$ holds, so we substitute
$(\beta-\alpha)_{+}=0$ and $\tilde{\rho}=(1+\half\epsilon_{2})^{-1}\sqrt{r/r^{\star}}$,
reparameterize $\alpha=\lambda\tilde{\beta}$ with $\lambda\ge1$,
and solve 
\begin{gather*}
\min_{\begin{subarray}{c}
\alpha\ge\tbeta\ge0,\\
\alpha^{2}+\tilde{\rho}^{2}\tilde{\beta}^{2}\le1
\end{subarray}}\frac{1-2\alpha\tilde{\beta}+\tilde{\beta}^{2}}{1-\tilde{\beta}^{2}}=\min_{\begin{subarray}{c}
\tbeta\ge0,\lambda\ge1,\\
(\lambda^{2}+\tilde{\rho}^{2})\tilde{\beta}^{2}\le1
\end{subarray}}\frac{1-2\lambda\tilde{\beta}^{2}+\tilde{\beta}^{2}}{1-\tilde{\beta}^{2}}=\min_{\lambda\ge1}\frac{(\lambda^{2}+\tilde{\rho}^{2})-2\lambda+1}{(\lambda^{2}+\tilde{\rho}^{2})-1}\\
=1-2\max_{\lambda\ge1}\frac{(\lambda-1)}{(\lambda^{2}-1)+\tilde{\rho}^{2}}=1-2\left(1+\min_{\lambda\ge1}\left\{ \lambda+\frac{\tilde{\rho}^{2}}{\lambda-1}\right\} \right)^{-1}=\frac{1}{1+\tilde{\rho}^{-1}}.
\end{gather*}
This yields $\delta^{\star}=(1+\tilde{\rho}^{-1})^{-1}=(1+(1+\frac{1}{2}\epsilon_{2})\sqrt{r^{\star}/r})^{-1}$
as claimed.
\end{proof}
The proof is completed by using \lemref{trivial1} to cover the two
trivial cases, and substituting \lemref{techical} into \lemref{nontriv1}
to cover the nontrivial cases.
\begin{theorem}
\label{thm:sym_detail}Fix $X\in\R^{n\times r}$ and $Z\in\R^{n\times r^{\star}}$
such that $r\ge r^{\star}$ and $XX^{T}\ne ZZ^{T}$. Suppose that
the RIP constant $0\le\delta<1$ and accuracy parameters $\epsilon_{0},\epsilon_{1},\epsilon_{2}\ge0$
satisfy the following:
\[
\epsilon_{0}\sqrt{r+r^{\star}}+\frac{1}{2}\epsilon_{1}-\left(\frac{1}{1+(1+\half\epsilon_{2})\sqrt{r^{\star}/r}}-\delta\right)\|XX^{T}-ZZ^{T}\|<0.
\]
Then, for every $f\in\cF_{Z}(\delta,\epsilon_{0})$, the universal
escape directions $D_{0},D_{1},\dots,D_{r^{\star}}$ defined in \defref{directions}
satisfy at least one of the following
\[
\frac{\inner{\nabla f(X)}{D_{0}}}{\|D_{0}X^{T}\|}<-\epsilon_{1}\quad\text{ and/or }\quad\min_{i\in\{1,\dots,r^{\star}\}}\frac{\inner{\nabla^{2}f(X)[D_{i}]}{D_{i}}}{\|D_{i}X^{T}\|^{2}}<-\epsilon_{2}.
\]
\end{theorem}

\begin{proof}
Recall that our goal is to demonstrate, for every $X,Z$ pair with
large error $E=XX^{T}-ZZ^{T}\ne0$, that there exists a choice of
$t_{1}\ge0,t_{2}\ge0$ to yield $\psi(t_{1},t_{2})<0$ in the dual
function defined in (\ref{eq:dual1}). First, for the trivial cases
of $X=0$ or $Z_{\perp}=(I-XX^{\dagger})Z=0$, we have either $\psi(1,0)<0$
or $\psi(0,1)<0$ via \lemref{trivial1}. Otherwise, for the nontrivial
case of $X\ne0$ and $Z_{\perp}\ne0$, we verify that
\begin{align*}
 & \max_{0\le\tau\le\alpha}\sqrt{1-\alpha^{2}}\sqrt{1-\tau^{2}}+[\alpha-(1+\delta+\half\epsilon_{2})\beta]\tau-\delta\\
\overset{\text{(a)}}{\ge} & \max_{0\le\tau\le\alpha}\sqrt{1-\alpha^{2}}\sqrt{1-\tau^{2}}+\alpha\tau-(1+\half\epsilon_{2})\beta\tau-[1+(1+\half\epsilon_{2})\beta\tau]\delta\\
= & \max_{0\le\tau\le\alpha}\left(\frac{\sqrt{1-\alpha^{2}}\sqrt{1-\tau^{2}}+\alpha\tau-(1+\half\epsilon_{2})\beta\tau}{1+(1+\half\epsilon_{2})\beta\tau}-\delta\right)[1+(1+\half\epsilon_{2})\beta\tau]\\
\overset{\text{(b)}}{\ge} & \left(\frac{1}{1+(1+\half\epsilon_{2})\sqrt{r^{\star}/r}}-\delta\right)>0.
\end{align*}
Step (a) is by substituting $\tau\beta\delta\le\tau(1+\half\epsilon_{2})\beta\delta$.
Step (b) follows by substituting \lemref{techical}; the first factor
in the product is nonnegative by hypothesis, and this allows us to
substitute $(1+\half\epsilon_{2})\beta\tau\ge0$ into the second factor.
Substituting the above into \lemref{nontriv1} yields a negative dual
function
\begin{multline*}
\min_{0\le\tau\le\alpha}\psi(\sqrt{\frac{1-\tau^{2}}{1-\alpha^{2}}},\frac{\tau}{\alpha})\le(\epsilon_{0}\sqrt{r+r^{\star}}+\half\epsilon_{1})\|E\|\\
-\|E\|^{2}\max_{0\le\tau\le\alpha}\left\{ \sqrt{1-\alpha^{2}}\sqrt{1-\tau^{2}}+[\alpha-(1+\delta+\half\epsilon_{2})\beta]\tau-\delta\right\} <0.
\end{multline*}
 
\end{proof}
\thmref{sym} immediately follows from \thmref{sym_detail}. Indeed,
if $X$ is instead an approximate second-order point that satisfies
\[
\inner{\nabla f(X)}D\ge-\epsilon_{1}\|DX^{T}\|,\quad\inner{\nabla^{2}f(X)[D]}D\ge-\epsilon_{2}\|DX^{T}\|^{2}\quad\text{for all }D,
\]
then none of the escape directions $D_{0},D_{1},\dots,D_{r^{\star}}$
will be able to escape $X$. By reversing the implications of \thmref{sym_detail},
we conclude that $X$ must recover the ground truth with a recovery
error of 
\[
\left(\epsilon_{0}\sqrt{r+r^{\star}}+\frac{1}{2}\epsilon_{1}\right)\left(\frac{1}{1+(1+\half\epsilon_{2})\sqrt{r^{\star}/r}}-\delta\right)_{+}^{-1}\ge\|XX^{T}-ZZ^{T}\|.
\]

We now convert \thmref{sym_detail} from the local norm into the Euclidean
norm. The key idea is to appeal to the classical observation for
the Burer--Monteiro factorization that if a second-order point is
also rank-deficient, then it is globally optimal. 
\begin{lemma}
\label{lem:rankdef1}Let $f\in\cF_{Z}(\delta,\epsilon_{0})$ and let
$X$ satisfy $\|\nabla f(X)\|\le\epsilon_{1}$ and $\nabla^{2}f(X)\succeq-\epsilon_{2}I$.
Then,
\begin{gather*}
(1-\delta)\|E\|^{2}\le\epsilon_{0}\sqrt{r+r^{\star}}\|E\|+\frac{\epsilon_{1}}{2}\|X\|+\frac{\epsilon_{2}}{2}\|Z\|^{2}+2(1+\delta)\lambda_{\min}(X^{T}X)\|Z\|^{2}
\end{gather*}
where $E=XX^{T}-ZZ^{T}$. 
\end{lemma}

\begin{proof}
\global\long\def\tD{\tilde{D}}%
Define $\tD_{0}=-\half X$ and $\tD_{i}=Ze_{i}v^{T}$, so that $E=-(X\tD_{0}^{T}+\tD_{0}X^{T})-\sum_{i}\tD_{i}\tD_{i}^{T}$.
Observe that
\begin{align*}
\|E\|_{\HH}^{2} & =\inner{-S}E+\inner{\HH(E)+S}{XX^{T}}-\inner{\HH(E)+S}{ZZ^{T}}\\
 & =-\inner SE-\inner{\nabla f(X)}{\tD_{0}}-\half\sum_{i=1}^{r^{\star}}\left[\inner{\nabla^{2}f(X)[\tD_{i}]}{\tD_{i}}-\|X\tD_{i}^{T}+\tD_{i}X^{T}\|_{\HH}^{2}\right]\\
 & \le\epsilon_{0}\|E\|_{\nuc}+\half\epsilon_{1}\|\tD_{0}\|+\half\sum_{i=1}^{r^{\star}}\left[\epsilon_{2}\|\tD_{i}\|^{2}+4(1+\delta)\|\tD_{i}X^{T}\|^{2}\right]\\
 & =\epsilon_{0}\|E\|_{\nuc}+\half\epsilon_{1}\|X\|+\half\|Z\|^{2}\left[\epsilon_{2}+4(1+\delta)\lambda_{\min}(X^{T}X)\right].
\end{align*}
Finally, it follows from $\rank(E)\le r+r^{\star}$ that $\|E\|_{\HH}^{2}\ge(1-\delta)\|E\|^{2}$
and $\|E\|_{\nuc}\le\sqrt{r+r^{\star}}\|E\|$. 
\end{proof}
Combining the two bounds in \lemref{rankdef1} and \thmref{sym_detail}
then yields the following bound.
\begin{corollary}
\label{cor:sym_detail}Let $f\in\cF_{Z}(\delta,\epsilon_{0})$ and
let $X$ satisfy $\|\nabla f(X)\|\le\epsilon_{1}$ and $\nabla^{2}f(X)\succeq-\epsilon_{2}$.
If $\Delta=\frac{1}{1+\sqrt{r^{\star}/r}}-\delta>0$, then
\begin{gather*}
\|XX^{T}-ZZ^{T}\|\le\frac{\epsilon_{0}\sqrt{r+r^{\star}}}{\Delta}+\sqrt{\frac{\epsilon_{1}\sqrt{\chi}R+\epsilon_{2}\chi R^{2}}{\Delta}}
\end{gather*}
where $R=\max\{\|X\|,\|Z\|\}$ and $\chi=(1+\delta)(1+\sqrt{r/r^{\star}})$.
\end{corollary}

\begin{proof}
Write $\lambda_{\min}\eqdef\lambda_{\min}(X^{T}X)$. If $\|E\|^{2}/\lambda_{\min}\ge2\chi R^{2}$,
then substituting $2(1+\delta)R^{2}\lambda_{\min}\le\frac{1}{1+\sqrt{r/r^{\star}}}\|E\|^{2}$
into \lemref{rankdef1} yields
\[
\left(1-\delta-\frac{1}{1+\sqrt{r/r^{\star}}}\right)\|E\|^{2}-\epsilon_{0}\sqrt{r+r^{\star}}\|E\|\le\frac{\epsilon_{1}}{2}R+\frac{\epsilon_{2}}{2}R^{2}.
\]
Otherwise, if $\|E\|^{2}/\lambda_{\min}\le2\chi R^{2}$, then substituting
$-\|DX^{T}\|^{2}\le-\lambda_{\min}\|D\|^{2}$ and the following
\begin{align*}
\frac{1}{1+(1+\half\frac{\epsilon_{2}}{\lambda_{\min}})\sqrt{r^{\star}/r}} & =1-\frac{(1+\half\frac{\epsilon_{2}}{\lambda_{\min}})\sqrt{r^{\star}/r}}{1+(1+\half\frac{\epsilon_{2}}{\lambda_{\min}})\sqrt{r^{\star}/r}}\\
 & \ge1-\frac{\sqrt{r^{\star}/r}}{1+\sqrt{r^{\star}/r}}-\frac{\half\frac{\epsilon_{2}}{\lambda_{\min}}\sqrt{r^{\star}/r}}{\sqrt{r/r^{\star}}}=1-\frac{1}{1+\sqrt{r/r^{\star}}}-\half\frac{\epsilon_{2}}{\lambda_{\min}}
\end{align*}
into \thmref{sym} yields the following, after substituting $\|E\|^{2}/\lambda_{\min}\le2\chi R^{2}$
\begin{align*}
\left(1-\delta-\frac{1}{1+\sqrt{r/r^{\star}}}\right)\|E\|^{2}-\epsilon_{0}\sqrt{r+r^{\star}}\|E\| & \le\frac{\epsilon_{1}\sqrt{2\chi}R}{2}+\frac{\epsilon_{2}2\chi R^{2}}{2}.
\end{align*}
Finally, solving $a\|E\|^{2}-b\|E\|-c\le0$ yields $a\|E\|\le\frac{1}{2}(b+\sqrt{b^{2}+4ac})\le b+\sqrt{ac}$
and hence $\|E\|\le b/a+\sqrt{c/a}$.
\end{proof}
We conclude the proof of \corref{sym} by substituting the following
slight refinement of \cite[Lemma~1.1]{candes2011tight}, which explicitly
spells out the leading constants. Its proof is deferred to \appref{specnorm}.
\begin{lemma}
\label{lem:specnorm}Let $\AA:\R^{n_{1}\times n_{2}}\to\R^{m}$ satisfy
$(1,\delta)$-RIP for $\delta<1$. If $z\sim\mathcal{N}(0,\sigma^{2}I_{m})$,
then $\|\AA^{T}(z)\|_{\op}\le10\sigma\sqrt{n_{1}+n_{2}}$ holds with
probability $1-12^{-n_{1}-n_{2}+1}$.
\end{lemma}

\subsection{\protect\label{sec:asym}Asymmetric parameterization $UV^{T}$}

Let us now see how the same approach generalizes in a natural way
to the following class of asymmetric recovery problems with a balancing
regularizer
\[
\cG_{Z}(\delta,\epsilon_{0})=\left\{ g:\R^{n\times r}\to\R:\begin{array}{c}
g([U;V])=2\|\AA(UV^{T})-b\|^{2}+\half\|U^{T}U-V^{T}V\|^{2},\\
\|\AA^{T}(\AA(U_{\star}V_{\star}^{T})-b)\|_{\op}\le\epsilon_{0}\quad\text{for }Z=[U_{\star};V_{\star}],\\
\AA\in\ripd.
\end{array}\right\} 
\]
The key idea is to rewrite the asymmetric problem over $n_{1}\times n_{2}$
low-rank matrices into a regularized version of the symmetric problem
over $(n_{1}+n_{2})\times(n_{1}+n_{2})$ low-rank matrices
\begin{align*}
g([U;V]) & =2\|\AA(UV^{T})-b\|^{2}+\frac{1}{2}\|U^{T}U-V^{T}V\|^{2}\\
 & =\frac{1}{2}\|\BB(XX^{T})-2b\|^{2}+\frac{1}{2}\|X^{T}JX\|^{2}.
\end{align*}
Here, $J=\diag(I_{n_{1}},-I_{n_{2}})$, and the linear operator $\BB$
is explicitly defined as follows 
\[
[\BB(XX^{T})]_{i}=\inner{B_{i}}{XX^{T}}=\inner{\begin{bmatrix}0 & A_{i}\\
A_{i}^{T} & 0
\end{bmatrix}}{\begin{bmatrix}U\\
V
\end{bmatrix}\begin{bmatrix}U\\
V
\end{bmatrix}^{T}}=2\inner{A_{i}}{UV^{T}}=2[\AA(UV^{T})]_{i}
\]
for each $i=\{1,2,\dots m\}$. This regularizer is needed because
$\BB$ does not actually satisfy RIP on its own. Instead, the regularizer
term $J$ is needed up to ``prop up'' the null space of $\BB$. 

\begin{fact}
Define $\HH$ such that $\HH(E)=\BB^{T}\BB(E)+JEJ$ for all $E\in\S^{n}$.
If $\AA\in\rip(\delta,k)$, then $\left|\|M\|_{\HH}^{2}/\|M\|^{2}-1\right|\le2\delta$
holds for all $M\in\S^{n}$ such that $\rank(M)\le k$.
\end{fact}

\begin{proof}
We make the following partition
\[
E=\begin{bmatrix}E_{11} & E_{12}\\
E_{12}^{T} & E_{22}
\end{bmatrix},\quad JEJ=\begin{bmatrix}E_{11} & -E_{12}\\
-E_{12}^{T} & E_{22}
\end{bmatrix},\quad\BB(E)=2\AA(E_{12}),
\]
and observe that $\rank(E_{12})\le\rank(E)$, and hence
\begin{align*}
\|\BB(E)\|^{2}+\inner E{JEJ} & \ge4(1-\delta)\|E_{12}\|^{2}+\|E_{11}\|^{2}-2\|E_{12}\|^{2}+\|E_{22}\|^{2}\\
 & =\|E_{11}\|^{2}+2(1-2\delta)\|E_{12}\|^{2}+\|E_{22}\|^{2}\\
 & \ge\left(1-2\delta\right)\left(\|E_{11}\|^{2}+2\|E_{12}\|^{2}+\|E_{22}\|^{2}\right)=(1-2\delta)\|E\|^{2}.
\end{align*}
The proof of the upper-bound is identical.
\end{proof}
Moreover, the embedding does not affect the effect of noise.
\begin{fact}
For all $w\in\R^{m}$, we have $\|\BB^{T}(w)\|_{\op}=\|\AA^{T}(w)\|_{\op}$
. 
\end{fact}

\begin{proof}
We have $\|\BB^{T}(w)\|_{\op}=\max_{i}|\lambda_{i}[\BB^{T}(w)]|$
and $\|\AA^{T}(w)\|_{\op}=\max_{i}\sigma_{i}[\BB^{T}(w)]$ and the
following 
\[
\lambda_{i}[\BB^{T}(w)]=\lambda_{i}\left(\begin{bmatrix}0 & \AA^{T}(w)\\{}
[\AA^{T}(w)]^{T} & 0
\end{bmatrix}\right)=\pm\sigma_{i}[\AA^{T}(w)].
\]
\end{proof}
The remainder of the proof now closely parallels that of \thmref{sym_detail}.
For every balanced ground truth $Z=[U_{\star};V_{\star}]$ such that
$Z^{T}JZ=U_{\star}^{T}U_{\star}-V_{\star}^{T}V_{\star}=0$, and for
every spurious point $X=[U;V]$ with $XX^{T}\ne ZZ^{T}$, we will
identify a choice of rescaling $t_{1}\ge0,t_{2}\ge0$ so that $\psi(t_{1},t_{2})<0$
holds where
\begin{equation}
\psi(t_{1},t_{2})\eqdef\sup_{g\in\cG_{Z}(\delta,\epsilon_{0})}\left\{ \begin{array}{c}
t_{1}(\inner{\nabla g(X)}{D_{0}}+\epsilon_{1}\|D_{0}X^{T}\|)\\
+\half t_{2}\sum_{i=1}^{nr}(\inner{\nabla^{2}g(X)[D_{i}]}{D_{i}}+\epsilon_{2}\|D_{i}X^{T}\|^{2})
\end{array}\right\} .\label{eq:dual2}
\end{equation}
We again begin by stating the concrete expressions for the directional
derivatives of $g\in\cG_{Z}(\delta,\epsilon_{0})$ in terms of $g(X)=\frac{1}{2}\|\BB(XX^{T})-2b\|^{2}+\frac{1}{2}\|X^{T}JX\|^{2}$.
\begin{fact}
\label{fact:deriv2}For $g\in\cG_{Z}(\delta,\epsilon_{0})$, write
$S=\BB^{T}(\BB(ZZ^{T})-2b)$ and $\HH(M)\eqdef\BB^{T}\BB(M)+JMJ$.
Then, the directional derivatives of $g$ are written
\begin{align*}
\inner{\nabla g(X)}D & =\inner{\HH(E)+JZZ^{T}J+S}{DX^{T}+XD^{T}},\\
\half\inner{\nabla^{2}g(X)[D_{i}]}D & =\inner{\HH(E)+JZZ^{T}J+S}{DD^{T}}+\half\|DX^{T}+XD^{T}\|_{\HH}^{2}.
\end{align*}
Moreover, $\|S\|_{\op}\le2\epsilon_{0}$ and $\left|\|M\|_{\HH}^{2}/\|M\|^{2}-1\right|\le2\delta$
for all $E\in\S^{n}$ such that $\rank(E)\le r+r^{\star}$.
\end{fact}

The following two lemmas are analogs of \lemref{trivial1} and \lemref{nontriv1}.
Their proofs are essentially verbatim, except that \lemref{slacklb}
and $(JZZ^{T}J)Z=JZ(Z^{T}JZ)=0$ are used to lower-bound $\inner{JZZ^{T}J}{t_{1}E_{\T}+t_{2}E_{\N}}\ge0$
for $t_{1}\ge t_{2}\ge0$. 
\begin{lemma}
\label{lem:trivial2}Fix $X\in\R^{n\times r}$ and $Z\in\R^{n\times r^{\star}}$
such that $r\ge r^{\star},$ $E=XX^{T}-ZZ^{T}\ne0,$ and $Z^{T}JZ=0$.
Define the dual function $\psi(t_{1},t_{2})$ as in (\ref{eq:dual2}).
If $X\ne0$ and $XX^{\dagger}Z=Z$, then 
\[
\psi(1,0)\le(2\sqrt{r+r^{\star}}\epsilon_{0}+\half\epsilon_{1})\|E\|-(1-2\delta)\|E\|^{2}.
\]
If instead $X=0$, and hence $XX^{\dagger}Z=0$, then 
\[
\psi(0,1)\le2\sqrt{r+r^{\star}}\epsilon_{0}\|E\|-(1-2\delta)\|E\|^{2}.
\]
\end{lemma}

\begin{lemma}
\label{lem:nontrivial2}Fix $X\in\R^{n\times r}$ and $Z\in\R^{n\times r^{\star}}$
such that $r\ge r^{\star},$ $X\ne0,$ $E=XX^{T}-ZZ^{T}\ne0,$ $Z_{\perp}=(I-XX^{\dagger})\ne0,$
and $Z^{T}JZ=0$. Define the parameters $\alpha,\beta$ as follows
\[
\alpha=\frac{\|Z_{\perp}Z_{\perp}^{T}\|}{\|E\|},\qquad\beta=\frac{\lambda_{\min}(X^{T}X)}{\|E\|}\frac{\tr(Z_{\perp}Z_{\perp}^{T})}{\|Z_{\perp}Z_{\perp}^{T}\|}.
\]
Then, the dual function $\psi(t_{1},t_{2})$ in (\ref{eq:dual2})
is upper-bounded
\begin{align*}
\psi(\sqrt{\frac{1-\tau^{2}}{1-\alpha^{2}}},\frac{\tau}{\alpha}) & \le(2\epsilon_{0}\sqrt{r+r^{\star}}+\half\epsilon_{1})\|E\|\\
 & \qquad-\left[\sqrt{1-\alpha^{2}}\sqrt{1-\tau^{2}}+[\alpha-(1+2\delta+\half\epsilon_{2})\beta]\tau-2\delta\right]\|E\|^{2}.
\end{align*}
\end{lemma}

We now repeat the proof of \thmref{sym_detail} near-verbatim. The
only superficial differences are: 1) the RIP constant $\delta$ and
the noise level $\epsilon_{0}$ are effectively doubled; 2) the restriction
of the maximization over $\tau$ to $0\le\tau\le\alpha$ in \lemref{techical}
is explicitly required by \lemref{nontrivial2}; 3) the additional
requirement for $Z^{T}JZ=0$ can simply be relaxed when optimizing
over $X$ and $Z$ in \lemref{techical}. 
\begin{theorem}
\label{thm:nonsym_detail}Fix $X\in\R^{n\times r}$ and $Z\in\R^{n\times r^{\star}}$
such that $r\ge r^{\star}$, $XX^{T}\ne ZZ^{T},$ and $Z^{T}JZ=0$.
Suppose that the RIP constant $0\le\delta<1$ and accuracy parameters
$\epsilon_{0},\epsilon_{1},\epsilon_{2}\ge0$ satisfy the following:
\[
\frac{1}{2}\epsilon_{1}+2\epsilon_{0}\sqrt{r+r^{\star}}-\left(\frac{1}{1+(1+\half\epsilon_{2})\sqrt{r^{\star}/r}}-2\delta\right)\|XX^{T}-ZZ^{T}\|<0.
\]
Then, for every $g\in\cG_{Z}(\delta,\epsilon_{0})$, the escape directions
$D_{0},D_{1},\dots,D_{r^{\star}}$ defined in \defref{directions}
satisfy one of the following
\[
\frac{\inner{\nabla g(X)}{D_{0}}}{\|D_{0}X^{T}\|}<-\epsilon_{1}\quad\text{ or }\quad\min_{i\in\{1,\dots,r^{\star}\}}\frac{\inner{\nabla^{2}g(X)[D_{i}]}{D_{i}}}{\|D_{i}X^{T}\|^{2}}<-\epsilon_{2}.
\]
\end{theorem}

Finally, converting \thmref{nonsym_detail} from the local norm into
the Euclidean norm requires the following lemma, which is the natural
analog for \lemref{rankdef1}. 
\begin{lemma}
\label{lem:rankdef2}Let $g\in\cG_{Z}(\delta,\epsilon_{0})$ where
$Z^{T}JZ=0$, and let $X$ satisfy $\|\nabla g(X)\|\le\epsilon_{1}$
and $\nabla^{2}g(X)\succeq-\epsilon_{2}I$. Then,
\begin{gather*}
(1-2\delta)\|E\|^{2}\le2\epsilon_{0}\sqrt{r+r^{\star}}\|E\|+\frac{\epsilon_{1}}{2}\|X\|+\frac{\epsilon_{2}}{2}\|Z\|^{2}+2(1+2\delta)\lambda_{\min}(X^{T}X)\|Z\|^{2}.
\end{gather*}
where $E=XX^{T}-ZZ^{T}$. 
\end{lemma}

The proof of the following then follows by repeating the proof of
\corref{sym_detail} verbatim, but using \lemref{rankdef2} and \thmref{nonsym_detail}.
\begin{corollary}
\label{cor:nonsym_detail}Let $g\in\cG_{Z}(\delta,\epsilon_{0})$
where $Z^{T}JZ=0$, and let $X$ satisfy $\|\nabla g(X)\|\le\epsilon_{1}$
and $\nabla^{2}g(X)\succeq-\epsilon_{2}$. Then,
\begin{gather*}
\|XX^{T}-ZZ^{T}\|\le\frac{2\epsilon_{0}\sqrt{r+r^{\star}}}{\Delta}+\sqrt{\frac{\epsilon_{1}\sqrt{\chi}R+\epsilon_{2}\chi R^{2}}{\Delta}}
\end{gather*}
where $\Delta=\frac{1}{1+\sqrt{r^{\star}/r}}-2\delta$ and $R=\max\{\|X\|,\|Z\|\}$
and $\chi=(1+2\delta)(1+\sqrt{r/r^{\star}})$.
\end{corollary}

We conclude the proof of \corref{nonsym} by substituting \lemref{specnorm},
which says that $\|\AA^{T}(w)\|_{\op}\le10\sigma\sqrt{n_{1}+n_{2}}$
holds with high probability when $w\sim\mathcal{N}(0,\sigma^{2}I_{m})$.

\section{\protect\label{sec:counter}Proof of necessity by counterexamples}

The counterexample we use to establish the sharpness of \thmref{sym}
reads as follows. 
\begin{example}[Symmetric parameterization]
\label{exa:sym}Let $[Q_{1},Q_{2}]$ have orthonormal columns with
$Q_{1}\in\R^{n\times r}$ and $Q_{2}\in\R^{n\times r^{\star}}$ satisfying
$1\le r^{\star}\le r<n$. For $\varepsilon\ge0$, let
\[
\delta=\frac{1}{1+(1+\varepsilon)\sqrt{r^{\star}/r}},\qquad M^{\star}=(1+\varepsilon+\sqrt{r/r^{\star}})Q_{2}Q_{2}^{T},
\]
and implicitly define $\AA:\R^{n\times n}\to\R^{n\times n}$ to satisfy
the following 
\begin{gather*}
\|\AA(E)\|^{2}=(1+\delta)\|E\|^{2}-2\delta\inner GE^{2}\quad\text{for all }E\in\R^{n\times n}\\
\text{where }G=\frac{1}{\sqrt{2r}}Q_{1}Q_{1}^{T}-\frac{1}{\sqrt{2r^{\star}}}Q_{2}Q_{2}^{T}.
\end{gather*}
It follows from $\|G\|=1$ that $\AA\in\rip(\delta,n)\subseteq\rip(\delta,k)$
for any $k$. However, the function $f(X)=\frac{1}{2}\|\AA(XX^{T}-M^{\star})\|^{2}$
has a spurious approximate second-order point at the point $X_{0}=Q_{1}$:
\[
\nabla f(X_{0})=0,\quad\nabla^{2}f(X_{0})\succeq-2(1+\delta)\varepsilon I,\quad\|X_{0}X_{0}^{T}-M^{\star}\|>\|M^{\star}\|.
\]
\end{example}

To verify the first- and second-order conditions in \exaref{sym},
we require the following lemma, which is a modified version of \cite[Lemma~6.1]{zhang2022improved}. 
\begin{lemma}
\label{lem:counter}Under the same setting as \exaref{sym}, let $\H\in\S^{n^{2}\times n^{2}}$
satisfy the following $rr^{\star}+2$ eigenvalue equations
\begin{gather*}
\H\vect\left(Q_{1}Q_{1}^{T}+\sqrt{r/r^{\star}}Q_{2}Q_{2}\right)=\vect\left(Q_{1}Q_{1}^{T}+\sqrt{r/r^{\star}}Q_{2}Q_{2}\right),\\
\H\vect\left(Q_{1}Q_{1}^{T}-\sqrt{r/r^{\star}}Q_{2}Q_{2}\right)=\frac{1-\delta}{1+\delta}\cdot\vect\left(Q_{1}Q_{1}^{T}-\sqrt{r/r^{\star}}Q_{2}Q_{2}\right),\\
\H\vect\left(Q_{1}VQ_{2}^{T}+Q_{2}V^{T}Q_{1}^{T}\right)=\vect\left(Q_{1}VQ_{2}^{T}+Q_{2}V^{T}Q_{1}^{T}\right)\qquad\text{for all }V\in\R^{r\times r^{\star}}.
\end{gather*}
Then, the function $h(X)\eqdef\frac{1}{2}\|\H^{1/2}\vect(XX^{T}-M^{\star})\|^{2}$
satisfies $\nabla h(Q_{1})=0$ and $\nabla^{2}h(Q_{1})\succeq-2\varepsilon I$.
\end{lemma}

\begin{proof}
Write $\rho\equiv\sqrt{r/r^{\star}}$ and $\tilde{\rho}=(1+\varepsilon)^{-1}\rho$,
and take $[Q_{1},Q_{2}]=I_{r+r^{\star}}$ without loss of generality.
Decomposing $\e=\vect(XX^{T}-M^{\star})$ into eigenvectors and applying
$\H$ yields
\begin{gather*}
\mat(\e)=\begin{bmatrix}I_{r} & 0\\
0 & -(1+\varepsilon+\rho)I_{r^{\star}}
\end{bmatrix}=-\frac{1}{2\tilde{\rho}}\begin{bmatrix}I_{r} & 0\\
0 & \rho I_{r^{\star}}
\end{bmatrix}+\frac{1+2\tilde{\rho}}{2\tilde{\rho}}\begin{bmatrix}I_{r} & 0\\
0 & -\rho I_{r^{\star}}
\end{bmatrix},\\
\mat(\H\e)=-\frac{1}{2\tilde{\rho}}\begin{bmatrix}I_{r} & 0\\
0 & \rho I_{r^{\star}}
\end{bmatrix}+\frac{1}{2\tilde{\rho}}\begin{bmatrix}I_{r} & 0\\
0 & -\rho I_{r^{\star}}
\end{bmatrix}=\begin{bmatrix}0_{r} & 0\\
0 & -(1+\varepsilon)I_{r^{\star}}
\end{bmatrix},
\end{gather*}
where we used $(1-\delta)/(1+\delta)=(1+2\tilde{\rho})^{-1}$. Now,
to verify that $\nabla h(Q_{1})=0$, one can check for arbitrary $V=[V_{1};V_{2}]$
with $V_{1}\in\R^{r\times r}$ and $V_{2}\in\R^{r^{\star}\times r}$
that
\begin{align*}
\inner{\nabla h(Q_{1})}V & =\inner{\mat(\H\e)}{Q_{1}V^{T}+VQ_{1}^{T}}\\
 & =\left\langle \begin{bmatrix}0_{r} & 0\\
0 & -(1+\varepsilon)I_{r^{\star}}
\end{bmatrix},\begin{bmatrix}V_{1}+V_{1}^{T} & V_{2}^{T}\\
V_{2} & 0_{r^{\star}}
\end{bmatrix}\right\rangle =0.
\end{align*}
Similarly, to verify that $\nabla^{2}h(Q_{1})\succeq-2\varepsilon I$,
one can check that
\begin{align*}
\inner{\nabla^{2}h(Q_{1})[V]}V & =2\inner{\mat(\H\e)}{VV^{T}}+\|\H^{1/2}\vect(Q_{1}V^{T}+VQ_{1}^{T})\|^{2}\\
 & \ge-2(1+\varepsilon)\|V_{2}\|^{2}+2\|V_{2}\|^{2}=-2\varepsilon\|V_{2}\|^{2}\ge-2\varepsilon\|V\|^{2}
\end{align*}
where the second line follows from the first because
\begin{gather*}
\langle\mat(\H\e),VV^{T}\rangle=\left\langle \begin{bmatrix}0_{r} & 0\\
0 & -(1+\varepsilon)I_{r^{\star}}
\end{bmatrix},\begin{bmatrix}V_{1}V_{1}^{T} & V_{1}V_{2}^{T}\\
V_{2}V_{1}^{T} & V_{2}V_{2}^{T}
\end{bmatrix}\right\rangle =-(1+\varepsilon)\|V_{2}\|^{2},\\
\left\Vert \H^{1/2}\vect\left(\begin{bmatrix}V_{1}+V_{1}^{T} & V_{2}^{T}\\
V_{2} & 0_{r^{\star}}
\end{bmatrix}\right)\right\Vert ^{2}\ge\left\Vert \H^{1/2}\vect\left(\begin{bmatrix}0_{r} & V_{2}^{T}\\
V_{2} & 0_{r^{\star}}
\end{bmatrix}\right)\right\Vert ^{2}=2\|V_{2}\|^{2},
\end{gather*}
where we recall that $\vect\left(\begin{bmatrix}0_{r} & V_{2}^{T}\\
V_{2} & 0_{r^{\star}}
\end{bmatrix}\right)$ is an eigenvector of $\H$ by hypothesis. 
\end{proof}
It turns out that the same counterexample for the symmetric case also
extends to the asymmetric case. The following lemma outlines the two
critical properties that allow us to extend \exaref{sym} from the
symmetric parameterization to the asymmetric parameterization.
\begin{lemma}
\label{lem:sym2nonsym}Let $f(X)=\|\AA(XX^{T}-M^{\star})\|^{2}$ for
$M^{\star}\in\S^{n}$, and let $X_{0}\in\R^{n\times r}$ be an approximate
second-order point with $\nabla f(X_{0})=0$ and $\nabla^{2}f(X_{0})\succeq-\epsilon I$.
If the following two conditions hold
\begin{gather}
\AA^{T}\AA(M^{T})=[\AA^{T}\AA(M)]^{T}\quad\text{for all }M\in\R^{n\times n},\label{eq:nonsym_nec_1}\\
\inner{\AA(X_{0}X_{0}^{T}-M^{\star})}{\AA(\dot{Y}\dot{Y}^{T})}\le0\quad\text{for all }\dot{Y}\in\R^{n\times r},\label{eq:nonsym_nec_2}
\end{gather}
then the function $g([U;V])=\|\AA(UV^{T}-M^{\star})\|^{2}$ with the
same $M^{\star}$ has $U_{0}=V_{0}=X_{0}$ as an approximate second-order
point with $\nabla g([U_{0};V_{0}])=0$ and $\nabla^{2}g([U_{0};V_{0}])\succeq-\epsilon I$. 
\end{lemma}

\begin{proof}
Let $S=\AA^{T}\AA(E)$ where $E=X_{0}X_{0}^{T}-M^{\star}=U_{0}V_{0}^{T}-M^{\star}$.
This same matrix $S$ appears in the first two directional derivatives
of $f$ at the point $X_{0}$, 
\begin{gather*}
\inner{\nabla f(X_{0})}{\dot{X}_{1}}=\inner S{\dot{X}_{1}X_{0}^{T}+X_{0}\dot{X}_{1}^{T}}=\inner{(S+S^{T})X_{0}}{\dot{X}_{1}},\\
\inner{\nabla^{2}f(X_{0})[\dot{X}_{1}]}{\dot{X}_{2}}=\inner S{\dot{X}_{1}\dot{X}_{2}^{T}+\dot{X}_{2}\dot{X}_{1}^{T}}+\inner{\AA(\dot{X}_{1}X_{0}^{T}+X_{0}\dot{X}_{1}^{T})}{\AA(\dot{X}_{2}X_{0}^{T}+X_{0}\dot{X}_{2}^{T})},
\end{gather*}
and also in the first two directional derivatives of $g$ at the point
$[U_{0};V_{0}]$,
\begin{gather*}
\inner{\nabla g(\begin{bmatrix}U_{0}\\
V_{0}
\end{bmatrix})}{\begin{bmatrix}\dot{U}_{1}\\
\dot{V}_{1}
\end{bmatrix}}=\inner S{\dot{U}_{1}V_{0}^{T}+U_{0}\dot{V}_{1}^{T}}=\inner{SV_{0}}{\dot{U}_{1}}+\inner{S^{T}U_{0}}{\dot{V}_{1}},\\
\inner{\nabla^{2}g(\begin{bmatrix}U_{0}\\
V_{0}
\end{bmatrix})[\begin{bmatrix}\dot{U}_{1}\\
\dot{V}_{1}
\end{bmatrix}]}{\begin{bmatrix}\dot{U}_{2}\\
\dot{V}_{2}
\end{bmatrix}}=\inner S{\dot{U}_{2}\dot{V}_{1}^{T}+\dot{U}_{1}\dot{V}_{2}^{T}}+\inner{\AA(\dot{U}_{1}V_{0}^{T}+U_{0}\dot{V}_{1}^{T})}{\AA(\dot{U}_{2}V_{0}^{T}+U_{0}\dot{V}_{2}^{T})}.
\end{gather*}
Under (\ref{eq:nonsym_nec_1}), it follows from $E=E^{T}$ that $S=S^{T}$.
If $\nabla f(X_{0})=0$, then $(S+S^{T})X_{0}=0$, and hence $SU_{0}=S^{T}V_{0}=0$,
and $\nabla g([U_{0};V_{0}])=0$. Next, we observe under (\ref{eq:nonsym_nec_1})
that
\begin{align*}
\inner{\nabla^{2}g(\begin{bmatrix}X_{0}\\
X_{0}
\end{bmatrix})[\begin{bmatrix}\dot{X}\\
\dot{X}
\end{bmatrix}]}{\begin{bmatrix}\dot{Y}\\
-\dot{Y}
\end{bmatrix}} & =\langle S,\dot{Y}\dot{X}^{T}-\dot{X}\dot{Y}^{T}\rangle+\inner{\AA(\dot{X}X_{0}^{T}+X_{0}\dot{X}^{T})}{\AA(\dot{Y}X_{0}^{T}-X_{0}\dot{Y}^{T})}\\
 & =0+\inner{\AA^{T}\AA(\dot{X}X_{0}^{T}+X_{0}\dot{X}^{T})}{\dot{Y}X_{0}^{T}-X_{0}\dot{Y}^{T}}=0
\end{align*}
because the inner product between a symmetric matrix and a skew-symmetric
matrix is always zero. For any arbitrary $\dot{U},\dot{V}\in\R^{n\times r}$,
let $\dot{X}=\frac{1}{2}(\dot{U}+\dot{V})$ and $\dot{Y}=\frac{1}{2}(\dot{U}-\dot{V})$.
Substituting $U_{0}=V_{0}=X_{0},$ $\dot{U}=\dot{X}+\dot{Y},$ and
$\dot{V}=\dot{X}-\dot{Y}$ yields
\[
\inner{\nabla^{2}g(\begin{bmatrix}U_{0}\\
V_{0}
\end{bmatrix})[\begin{bmatrix}\dot{U}\\
\dot{V}
\end{bmatrix}]}{\begin{bmatrix}\dot{U}\\
\dot{V}
\end{bmatrix}}=\inner{\nabla^{2}g(\begin{bmatrix}X_{0}\\
X_{0}
\end{bmatrix})[\begin{bmatrix}\dot{X}\\
\dot{X}
\end{bmatrix}]}{\begin{bmatrix}\dot{X}\\
\dot{X}
\end{bmatrix}}+\inner{\nabla^{2}g(\begin{bmatrix}X_{0}\\
X_{0}
\end{bmatrix})[\begin{bmatrix}\dot{Y}\\
-\dot{Y}
\end{bmatrix}]}{\begin{bmatrix}\dot{Y}\\
-\dot{Y}
\end{bmatrix}}.
\]
If $\nabla^{2}f(X_{0})\succeq-\epsilon I$, then the first of these
two terms is lower-bounded
\begin{align*}
\inner{\nabla^{2}g(\begin{bmatrix}X_{0}\\
X_{0}
\end{bmatrix})[\begin{bmatrix}\dot{X}\\
\dot{X}
\end{bmatrix}]}{\begin{bmatrix}\dot{X}\\
\dot{X}
\end{bmatrix}} & =2\inner S{\dot{X}\dot{X}^{T}}+\|\AA(\dot{X}V_{0}^{T}+U_{0}\dot{X}^{T})\|^{2}\\
 & =\inner{\nabla^{2}f(X_{0})[\dot{X}]}{\dot{X}}\ge-\epsilon\|\dot{X}\|^{2}\ge-\epsilon\|[\dot{U};\dot{V}]\|^{2}.
\end{align*}
where we used $\|\dot{U}+\dot{V}\|\le\sqrt{2}\|[\dot{U};\dot{V}]\|$.
Under (\ref{eq:nonsym_nec_2}), the second of these two terms is also
nonnegative
\[
\inner{\nabla^{2}g(\begin{bmatrix}U_{0}\\
V_{0}
\end{bmatrix})[\begin{bmatrix}\dot{Y}\\
-\dot{Y}
\end{bmatrix}]}{\begin{bmatrix}\dot{Y}\\
-\dot{Y}
\end{bmatrix}}=-2\inner S{\dot{Y}\dot{Y}^{T}}+\|\AA(\dot{Y}V_{0}^{T}-U_{0}\dot{Y}^{T})\|^{2}\ge0.
\]
\end{proof}
In turn, verifying that \exaref{sym} satisfies the two conditions
in \lemref{sym2nonsym} yields the following, which establishes the
sharpness of \thmref{asym}. 
\begin{example}[Asymmetric parameterization]
\label{exa:asym}Under the same setting as \exaref{sym}, the following
function
\[
g_{t}([U;V])=2\|\AA(UV^{T}-M^{\star})\|^{2}+t\|U^{T}U-V^{T}V\|^{2}\quad\text{for }t\ge0,
\]
admits a spurious second-order point at $U_{0}=V_{0}=Q_{1}$ that
satisfies
\[
\nabla g_{t}([U_{0};V_{0}])=0,\quad\nabla^{2}g_{t}([U_{0};V_{0}])\succeq-2(1+\delta)\varepsilon\cdot I,\quad\|U_{0}V_{0}^{T}-M^{\star}\|>\|M^{\star}\|.
\]
\end{example}

\begin{proof}
First, $\AA^{T}\AA(M^{T})=[\AA^{T}\AA(M)]^{T}$ holds because any
$\AA$ implicitly defined in \exaref{sym} must satisfy $\AA^{T}\AA(M)=(1+\delta)M-2\delta G\inner GM$
for $G=G^{T}$. Next, $S=\AA^{T}\AA(X_{0}X_{0}^{T}-M^{\star})=-(1+\delta)(1+\varepsilon)Q_{2}Q_{2}^{T}$
by repeating the proof of \lemref{counter}. Hence, the two conditions
in \lemref{sym2nonsym} hold. Finally, we observe that the twice-differentiable
function $\varphi([U;V])=\|U^{T}U-V^{T}V\|^{2}$ attains its minimum
at $U_{0}=V_{0}$; it must therefore satisfy $\nabla\varphi([U_{0};V_{0}])=0$
and $\nabla^{2}\varphi([U_{0};V_{0}])\succeq0$.
\end{proof}
The counterexample stated above is unaffected by the balancing regularizer.
However, if the balancing regularizer is eliminated, then spurious
approximate local minima can appear irrespective of overparameterization
and the value of the RIP constant. The following is a slight rephrasing
of \exaref{bal}. 
\begin{example}[Necessity of balancing regularizer]
For $r$ satisfying $1\le r<n$, define $g:\R^{2n\times r}\to\R$
such that $g([U;V])=\|UV^{T}-e_{1}e_{1}^{T}\|^{2}.$ Then, $U_{0}=(1/\sqrt{\epsilon})[e_{2},e_{3},\dots,e_{r+1}]$
and $V_{0}=0$ is a spurious approximate second-order point with 
\[
\nabla g([U_{0};V_{0}])=0,\quad\nabla^{2}g([U_{0};V_{0}])\succeq-\epsilon I,\quad\|U_{0}V_{0}^{T}-e_{1}e_{1}^{T}\|=1.
\]
\end{example}

\begin{proof}
The fact that $\nabla g([U_{0};V_{0}])=0$ follows from $V_{0}=0$
and $e_{1}^{T}U_{0}=0$, and therefore
\begin{gather*}
\inner{\nabla g(\begin{bmatrix}U_{0}\\
V_{0}
\end{bmatrix})}{\begin{bmatrix}\dot{U}\\
\dot{V}
\end{bmatrix}}=\inner{-e_{1}e_{1}^{T}}{\dot{U}_{1}V_{0}^{T}+U_{0}\dot{V}_{1}^{T}}=-\inner{e_{1}e_{1}^{T}U_{0}}{\dot{V}_{1}}=0.
\end{gather*}
Similarly, it follows from $V_{0}=0$ and $U_{0}^{T}U_{0}=\frac{1}{\epsilon}I_{r}$
that
\begin{align*}
\inner{\nabla^{2}g(\begin{bmatrix}U_{0}\\
V_{0}
\end{bmatrix})[\begin{bmatrix}\dot{U}\\
\dot{V}
\end{bmatrix}]}{\begin{bmatrix}\dot{U}\\
\dot{V}
\end{bmatrix}} & =2\inner{-e_{1}e_{1}^{T}}{\dot{U}\dot{V}^{T}}+\|\dot{U}V_{0}^{T}+U_{0}\dot{V}^{T}\|^{2}\\
 & =2\inner{-e_{1}e_{1}^{T}}{\dot{U}\dot{V}^{T}}+\frac{1}{\epsilon}\|\dot{V}\|^{2}=\inner{\begin{bmatrix}0 & -e_{1}e_{1}^{T}\\
-e_{1}e_{1}^{T} & \frac{1}{\epsilon}I
\end{bmatrix}}{\begin{bmatrix}\dot{U}\\
\dot{V}
\end{bmatrix}\begin{bmatrix}\dot{U}\\
\dot{V}
\end{bmatrix}^{T}}.
\end{align*}
The fact that $\nabla^{2}g([U_{0};V_{0}])\succeq-\epsilon$ holds
because
\[
\lambda_{\min}\left(\begin{bmatrix}0 & -1\\
-1 & \frac{1}{\epsilon}
\end{bmatrix}\right)=\frac{1}{2\epsilon}\left(1-\sqrt{1+4\epsilon^{2}}\right)\ge\frac{1}{2\epsilon}\left(1-1-2\epsilon^{2}\right)=-\epsilon.
\]
\end{proof}

\section{\protect\label{sec:conclusion}Concluding remarks}

This paper presents a unified, simplified, and strengthened proof
technique to derive sharp guarantees for nonconvex recovery in the
overparameterized regime. Using this technique, we fully complete
the picture between approximate second-order points and minimax-optimal
recovery under a rank-$(r+r^{\star})$ RIP assumption. An important
future work is to extend our analysis to non-RIP settings, and to
the rank-$k$ RIP setting where $k=\Omega(r^{\star})$. Based on our
results for the asymmetric parameterization, the most likely path
forward is by introducing further regularization. For example, incoherence
regularizers are used to make sparse measurements like matrix completion~\cite{keshavan2010matrixa,sun2016guaranteed,ge2017nospurious}
behave like dense RIP measurements, and nuclear norm regularizers
are used to select for low-rank solutions for undersampled measurements~\cite{candes2011tight,mcrae2024low}.
In turn, the success of specific algorithms like gradient descent
on the unregularized problem can again be attributed to their implicit
regularization. Indeed, our own analysis of the asymmetric case $UV^{T}$
worked precisely by augmenting it with an explicit regularizer to
reduce it to symmetric case $XX^{T}$ in the RIP setting. 

\subsubsection*{Acknowledgments}

I thank Simon S. Du, Fei Lu, Xiong Wang, and Mauro Maggioni for discussions
that motivated me to study the asymmetric case. Additionally, I thank
Salar Fattahi, C\'edric Josz, Andrew McRae, Sabrina Zielinski, and
two reviewers for helpful feedback and insightful suggestions. Financial
support for this work was provided by NSF CAREER Award ECCS-2047462
and ONR Award N00014-24-1-2671.

\bibliographystyle{siam}
\bibliography{proof_half}

\appendix

\section{\protect\label{app:specnorm}Proof of \lemref{specnorm}}

Denote $S(n)=\{x\in\R^{n}:\|x\|=1\}$ as the usual sphere. Write $Z=\AA^{T}(z)$
and assume $\sigma=1$ without loss of generality. Let $N_{\epsilon}(n)$
denote an $\epsilon$-net of $S(n)$ with respect to the Euclidean
norm, and note that $|N_{\epsilon}(n)|\le(3/\epsilon)^{n}$~\cite[Equation~III.1]{candes2011tight}.
Then, for every $u\in S(n_{1})$ and $v\in S(n_{2})$, we have
\begin{align*}
\inner u{Zv} & =\inner{u-u_{0}}{Zv}+\inner{Z^{T}u_{0}}{v-v_{0}}+\inner{u_{0}}{Zv_{0}}\\
 & \le\|Z\|_{\op}\|u-u_{0}\|+\|Z^{T}\|_{\op}\|v-v_{0}\|+\inner{u_{0}}{Zv_{0}}
\end{align*}
for some $u_{0}\in N_{1/4}(n_{1})$ and $v_{0}\in N_{1/4}(n_{2})$
such that $\|u-u_{0}\|\le1/4$ and $\|v-v_{0}\|\le1/4$. Hence, 
\[
\|Z\|_{\op}=\sup_{u\in S(n_{1}),v\in S(n_{2})}\inner u{Zv}\le2\sup_{u\in N_{1/4}(n_{1}),v\in N_{1/4}(n_{2})}\inner u{Zv}.
\]
For any fixed $u_{0},v_{0}$, we observe that the following is a standard
Gaussian variable
\[
\inner{u_{0}}{Zv_{0}}=\inner{u_{0}v_{0}^{T}}{\AA^{T}(z)}=\inner{\AA(u_{0}v_{0}^{T})}z\sim\mathcal{N}(0,\|\AA(u_{0}v_{0}^{T})\|^{2})
\]
whose variance is bounded $\|\AA(u_{0}v_{0}^{T})\|^{2}\le(1+\delta)\|u_{0}v_{0}^{T}\|^{2}\le2$.
Applying the standard Gaussian tail bound $\mathbb{P}(|\inner{u_{0}}{Zv_{0}}|\ge t)\le2\exp(-t^{2}/4)$
and taking the trivial union bound over all points on $N_{1/4}$ yields
\[
\mathbb{P}(\sup_{u_{0},v_{0}\in S_{1/4}}\inner{u_{0}}{Zu_{0}}\ge t)\le2|N_{1/4}(n_{1})||N_{1/4}(n_{2})|\exp\left(-\frac{t^{2}}{4}\right)\le2\exp\left((n_{1}+n_{2})\log12-\frac{t^{2}}{4}\right).
\]
Finally, if we set $t=5\sqrt{n_{1}+n_{2}}$, then $t^{2}/4\ge5(n_{1}+n_{2})>2(n_{1}+n_{2})\log12$,
and we obtain 
\begin{align*}
\mathbb{P}\left(\|Z\|_{\op}\ge10\sqrt{n_{1}+n_{2}}\right) & \le\mathbb{P}\left(2\sup_{u_{0},v_{0}\in S_{1/4}}\inner{u_{0}}{Zu_{0}}\ge10\sqrt{n_{1}+n_{2}}\right)\le2\exp\left(-(n_{1}+n_{2})\log12\right).
\end{align*}

\end{document}